\documentclass[leqno]{amsart}
\usepackage{amsmath,amsthm,amsfonts,amssymb}
\usepackage[backref=page]{hyperref}
\hypersetup{colorlinks=true,citecolor=blue,linkcolor=blue,urlcolor=blue,pdfstartview=FitH,
pdfauthor=Valentin Blomer ,pdftitle=Bounding sup-norms of cusp forms of large level}

\theoremstyle{plain}
\newtheorem{satz}{Theorem}
\newtheorem{lemma}{Lemma}
\newtheorem{prop}{Proposition}
\newtheorem{kor}{Corollary}

\theoremstyle{definition}

\renewcommand{\geq}{\geqslant}
\renewcommand{\leq}{\leqslant}

\begin{document}

\title{Bounding sup-norms of cusp forms of large level}

\author{Valentin Blomer}
\author{Roman Holowinsky}
\address{Department of Mathematics, University of Toronto, 40 St. George Street, Toronto, Ontario, Canada, M5S 2E4} \email{vblomer@math.toronto.edu} \email{romanh@math.toronto.edu}
\thanks{First author supported by NSERC grant 311664-05 and a Sloan Research Fellowship. Second author supported by NSERC and the Fields Institute Toronto.}


\begin{abstract} Let $f$ be an $L^2$-normalized  weight zero Hecke-Maa{\ss} cusp form of  square-free level $N$, character $\chi$ and Laplacian eigenvalue $\lambda\geq 1/4$. It is shown  that $\| f \|_{\infty} \ll_{\lambda} N^{-1/37}$, from which the hybrid bound $\|f \|_{\infty} \ll \lambda^{1/4} (N\lambda)^{-\delta}$ (for some $\delta > 0$) is   derived. The first bound holds also for $f = y^{k/2}F$ where $F$ is a holomorphic cusp form of weight $k$ with the implied constant now depending on $k$. 
 \end{abstract}


\maketitle

\section{Introduction}

If $f$ is a Maa{\ss} cusp form on a  (compactified) modular curve $X_0(N) =  \Gamma_0(N)\backslash \mathbb{H}^{\ast}$, normalized with respect to the  inner product
\begin{equation}\label{inner}
  \langle f_1, f_2 \rangle = \int_{X_0(N)} f_1(z) \bar{f}_2(z) \frac{dx dy}{y^2},
\end{equation}
it is interesting and often useful to derive pointwise bounds for $f$.  Such bounds play  an important role for instance in connection with subconvexity of $L$-functions \cite[Proposition 4]{HM}  and  in connection with the mass equidistribution conjecture \cite[Appendix]{Ru}. In the latter example, bounds for the $\Gamma_0(N)$-invariant function $f(z) = y^{k/2}|F(z)|$ are obtained for a holomorphic cusp form $F \in \mathcal{S}_k(N, \chi)$ of weight $k$. 

The problem of bounding cusp forms can be seen in the general framework of comparing two norms on a vector space. Of course, in the case of the finite-dimensional vector space $\mathcal{S}_k(N, \chi)$, the $L^2$-norm and the $L^{\infty}$-norm must be equivalent, but interestingly, the finite-dimensionality of $\mathcal{S}_k(N, \chi)$ follows conversely from the equivalence of these two norms \cite[Theorem 6.1]{Fr}. Another perspective comes from bounding periods of automorphic forms (see e.g.\  \cite{MV, Sa1}), and hence our work is at least implicitly related to the subconvexity problem for automorphic $L$-functions. We will come back to this point later.  In our case, of course, the ``period" is somewhat degenerate, namely a single point.

In the case $k=2$, the average
\begin{displaymath}
  \sup_{z \in \mathbb{H}} \sum_{F \in \mathcal{B}_2(N, 1)} y^2|F(z)|^2
\end{displaymath}
over an orthonormal basis $\mathcal{B}_2(N, 1)$ of $\mathcal{S}_2(N, 1)$ can be interpreted as the ratio of the Arakelov metric (coming from embedding $X_0(N)$ into its Jacobian) and the hyperbolic metric (coming from the Petersson inner product) on $X_0(N)$. Following \cite{AU, MU}, Jorgenson and Kramer \cite{JK} derived in a more general context the best-possible bound 
\begin{equation}\label{average}
  \sup_{z \in \mathbb{H}} \sum_{F \in \mathcal{B}_2(N, 1)} y^2|F(z)|^2 = O(1)
\end{equation}
as $N \rightarrow \infty$. Their line of investigation is continued in the important work \cite{JK2}. 

In order to treat the holomorphic and non-holomorphic case simultaneously, we define the archimedean parameter of a  cusp form $f$   as
\begin{equation}\label{arch}
  t = t_f = \left\{ \begin{array}{ll} 
    \frac{(k-1)}{2}, & f  = y^{k/2}F \text{ for a holomorphic form $F$ of weight } k,\\
    \sqrt{\lambda_f-\frac{1}{4}}, & f \text{ a Maa{\ss} form satisfying } \Delta f = \lambda_f f \text{ for } \Delta = -y^2(\partial_x^2 + \partial_y^2).
    \end{array}\right.
\end{equation}
The choice of the square root is irrelevant.   We call a Maa{\ss} form $f$
non-exceptional if $\lambda_f \geq 1/4$, i.e. $t_f \in \mathbb{R}$.  Pointwise bounds for $L^2$-normalized cusp forms can be  in terms of the archimedean parameter $t$, the level $N$, or both. The first breakthrough in the $t$-aspect was achieved by Iwaniec and Sarnak \cite{IS} who proved
\begin{equation}\label{IwSa}
  \| f \|_{\infty} \ll_{\varepsilon} (1+ |t_f|)^{5/12+\varepsilon} \| f \|_2.
\end{equation}
for a weight zero Hecke-Maa{\ss} cusp form of level 1, and any $\varepsilon > 0$. The ``trivial" bound in this context is $(1+|t_f|)^{1/2}$. One might conjecture that $\| f \|_{\infty} \ll_{\varepsilon} (1+|t_f|)^{\varepsilon} \|f\|_2$ holds for all $\varepsilon > 0$, but if at all, this can only hold for compact Riemann surfaces: high in the cusp one has \cite[p.\ 178]{Iw1}
\begin{displaymath}
  f\left(\frac{1}{4} + \frac{it_f}{2\pi}\right) = (1+|t_f|)^{1/6+o(1)} \| f \|_2.
\end{displaymath}
For holomorphic forms the situation is easier, because here the upper bound that comes from the Fourier expansion is essentially sharp. More precisely, it was recently shown in \cite{Xia} that
\begin{displaymath}
 k^{1/4-\varepsilon}\|f\|_2\ll_{\varepsilon}  \| f \|_{\infty} \ll_{\varepsilon} k^{1/4+\varepsilon}\|f\|_2
\end{displaymath}
for $f=y^{k/2}F$ where $F$ is a holomorphic Hecke cusp form of level 1 and weight $k$.  
For a sup-norm bound in a somewhat different context see also \cite{Ku}.\\
 
In the level aspect, nothing beyond the average bound \eqref{average} is known.  For an individual Hecke cusp form this  recovers only the trivial bound that comes from the Fourier expansion, see Lemma 7 below.  On the other hand, since the volume of $X_0(N)$ is about $N$, one might (optimistically) conjecture 
\begin{equation}\label{conjec}
  \| f \|_{\infty} \ll_{\varepsilon} N^{-1/2+\varepsilon} \|f\|_2.
\end{equation}
for fixed archimedean parameter $t$ and a Hecke eigenform $f$. If this was true, it is not only best possible, but would also be very strong in the sense  that it implied as a special case the Lindel\"of hypothesis for automorphic $L$-functions in the level aspect, since $f(i/\sqrt{N}) \approx L(1/2, f) N^{-1/2}  \| f \|_2$.  This shows, on the other hand, that in order to derive some subconvex bound for $L(1/2, f)$ in the level aspect by this method, one would need already a relatively strong pointwise bound $\|f\|_{\infty} \ll N^{-1/4-\delta} \|f\|_2$ for some $\delta > 0$. In a different direction, \eqref{conjec} would imply the best possible bounds 
\begin{equation}\label{lp}
  \| f \|_p \ll_{\varepsilon, p} N^{\frac{1}{p} - \frac{1}{2} + \varepsilon}, \quad p \geq 1.
\end{equation}
for the $L^p$-norm of an $L^2$-normalized Hecke cusp form $f$. 

Certainly \eqref{conjec} is true if $f$ is an oldform that comes from a level 1 form. It should be observed, however, that the requirement that $f$ is a Hecke eigenform cannot be dropped, as can be seen, for example, by looking at (holomorphic) Poincar\'e series
\begin{displaymath}
  P_m(z) :=   \sum_{n=1}^{\infty} \left(\delta_{nm} + 2\pi i^{-k} \sum_{N\mid c} \frac{1}{c}S_{\chi}(m, n; c) J_{k-1}\left(\frac{4\pi\sqrt{nm}}{c}\right)\right)\left(\frac{n}{m}\right)^{(k-1)/2}e(nz) \in \mathcal{S}_k(N, \chi).
\end{displaymath}
If we regard $m\geq 1$ and $k\geq 4$ as fixed and $N$ large, then $P_m(i) \asymp 1$, while it is easily seen that $\| y^{k/2} P_m \|_2 \asymp 1$, cf.\ e.g.\ \cite[(3.24)]{Iw2}.\\
 
The aim of this article is to establish the first non-trivial bound of the $L^{\infty}$-norm of an individual $L^2$-normalized holomorphic or non-holomorphic cuspidal Hecke eigenform of large square-free level $N$.  For notational simplicity let us write
\begin{displaymath}
  t^{\ast} := 1+|t_f|
\end{displaymath}
with $t_f$ as in \eqref{arch}. 

 \begin{satz} Let $N \geq 1$ be a square-free integer, $k \geq 2$ an even integer and $\chi$  an even Dirichlet character modulo $N$. Let $f$ be an $L^2$-normalized non-exceptional weight zero Hecke-Maa{\ss} cusp form of level $N$ and character $\chi$, or let $f=y^{k/2}F$ where $F \in \mathcal{S}_k(N, \chi)$ is an $L^2$-normalized Hecke eigenform. Then
 \begin{displaymath}
   \| f \|_{\infty} \ll (t^{\ast})^{\frac{11}{2}} N^{-\frac{1}{37}}.
 \end{displaymath}
 The implied constant is absolute.
 \end{satz}

As a first step towards \eqref{lp}, this implies
 
\begin{kor} Under the same assumptions as in Theorem 1, one has
\begin{displaymath}
 \| f \|_{p} \ll_{t} N^{-\frac{p-2}{37p}} 
\end{displaymath}
whenever $p \geq 2$. In particular, $\| f \|_4 \ll_t N^{-1/74}$. 
\end{kor}

As far as we know, no non-trivial bound for any $L^p$-norm, $p> 2$, in the $N$-aspect has ever been obtained. In the eigenvalue aspect, a very strong bound for the $L^4$-norm of level 1 Hecke-Maa{\ss} cusp forms has been proved by Sarnak and Watson, see  \cite[p.\ 460]{Sa} for more details and \cite{Sp} for the Eisenstein case. 
 
The method of the proof of Theorem 1 works also for holomorphic and non-holomorphic cusp forms of any weight with $(t^{\ast})^{11/2}$ replaced by $( |\lambda| + |k|)^A$. While we tried to be efficient in the $N$-aspect, we made little attempt to optimize the $t^{\ast}$-aspect. It is instructive to see that a nontrivial sup-norm bound is related to solving a very concrete diophantine problem (cf.\   Sections 7 and 8), and the focus here lies on a transparent treatment of the underlying number theoretical features rather than a tedious game with Bessel functions and oscillatory integrals. Our aim was, however,  to provide \emph{some} explicit $t^{\ast}$-exponent in Theorem 1 with the following application in mind. Analyzing the proof of \eqref{IwSa}, it is not hard to see that it is essentially  uniform in $N$: we will prove in Proposition 1 in Section 10 that
\begin{equation}\label{lastsection}
  \|f\|_{\infty} \ll_{\varepsilon} N^{\varepsilon}  (t^{\ast})^{5/12+\varepsilon} 
\end{equation}
for an $L^2$-normalized Hecke-Maa{\ss} cusp form of square-free level $N$. 
 Taking this for granted for the moment,  one obtains the following hybrid bound:

\begin{satz} If $f$ is a Maa{\ss} form as in Theorem 1, one has
\begin{displaymath}
  \| f \|_{\infty} \ll (t^{\ast})^{1/2}  (t^{\ast}N)^{-\frac{1}{2300}}
\end{displaymath}
with an absolute implied constant.
\end{satz}
 
Indeed, one has $\min((t^{\ast})^{5}N^{-\frac{1}{37}}, (t^{\ast})^{-\frac{1}{12}}) \leq ((t^{\ast})^{5}N^{-\frac{1}{37}})^{\frac{37}{2269}} ((t^{\ast})^{-\frac{1}{12}})^{\frac{2232}{2269}}  \leq (t^{\ast}N)^{-\frac{1}{2269}}$. \\

We briefly discuss the plan of attack. By Atkin-Lehner theory we can restrict $z$ to a Siegel set with $y \gg 1/N$ (here we use the fact that $N$ is square-free). Thus we will have to detect cancellation in sums of the type $\sum_{n \sim N} \lambda_f(n) e(nx)$ where we have no control on $x$, and little control on the Hecke eigenvalues $\lambda_f(n)$, since the level of $f$ is large. Generalizing \cite{MU}, we consider an amplified second moment. It is clear a priori that this setup alone cannot guarantee success, since in general the error term of the un-amplified second moment has about the same size as the main term. However, a careful analysis shows that for $x$ off rational points with small denominator,  this method does give a saving.   If $x$ is well-approximable, a Voronoi-type argument gives a non-trivial result.  The dichotomy of well and badly-approximable $x$ seems to be new in this context and does not appear in the work of Iwaniec-Sarnak \cite{IS}. It may  remind of the circle method, and we shall refer to the set of well-approximable $x$ as the major arcs.  

The paper is organized as follows: In Section 2-4 we recall some well-known material that we will need for the proof of Theorem 1. Section 5 gives a first bound for $f(x+iy)$ that will turn out to be useful for well-approximable $x$. Sections 6-9 implement the amplification method and give a second bound that is strong for badly-approximable $x$. We finish the proof of Theorem 1 in Section 9, and we prove \eqref{lastsection} in Section 10 which implies immediately Theorem 2.  Finally we remind the reader of the usual convention that the value of $\varepsilon$ may change during a calculation. \\

\textbf{Acknowledgement:} The authors thank Antal Balog, Matthew Young, Gergely Harcos, Philippe Michel, Jean Bourgain, Peter Sarnak and the referees for important and useful comments on various aspects of this work.

\section{Automorphic Forms and Fourier expansion}

Let $N\geq 1$ be a square-free integer, $\chi$  an even character to modulus
$N$; let $k\geq 2$ be an even integer. We denote by
  $\mathcal{L}^2(N,\chi)$ and
$\mathcal{L}^{2}_{0}(N,\chi)\subset \mathcal{L}^2(N,\chi)$, respectively, the spaces of weight zero Maa\ss\ forms, and of weight zero Maa\ss\ cusp forms with respect to the
congruence subgroup 
\begin{displaymath}
 \Gamma_{0}(N) = \left\{\left(\begin{smallmatrix}a & b \\ c & d \end{smallmatrix}\right)\in SL_2(\Bbb{Z}) : N \mid c  \right\}
\end{displaymath} 
  and with nebentypus $\chi$. In slight contrast to the usual notation, we shall write $\mathcal{S}_k(N, \chi)$ for the set of all functions $f(z) = y^{k/2}F(z)$ where $F$ is a holomorphic cusp form of weight $k$. Both $\mathcal{L}^2_0(N, \chi)$ and $\mathcal{S}_k(N, \chi)$ are Hilbert spaces with respect to the inner product \eqref{inner}. 
We can choose bases $\mathcal{B}_k(N, \chi)$ and $\mathcal{B}(N, \chi)$ of $\mathcal{S}_k(N, \chi)$ and $\mathcal{L}^2_0(N, \chi)$ that consist of $L^2$-normalized eigenforms for all Hecke operators $T_n$ with $(n, N) = 1$. For the proof of Theorem 1 we can assume that $f$ is a newform.  
 We define the archimedean parameter of a cusp form $f$ as in \eqref{arch} 
and write for notational simplicity
\begin{displaymath}
   t^{\ast} := 1+|t|
\end{displaymath}
as above. It is known \cite{KS} that
  $|\Im t | \leq \theta := \frac{7}{64}$. 
We write $\mathcal{B}(t, N, \chi) \subseteq \mathcal{B}(N, \chi)$ for the set of those $f \in \mathcal{B}(N, \chi)$ with archimedean parameter $t$. 
The orthogonal complement to $\mathcal{L}^2_{0}(N,\chi)$ in
$\mathcal{L}^2(N,\chi)$ is the Eisenstein spectrum $\mathcal{E}(N,\chi)$ (plus
possibly the space of constant functions if $\chi$ is trivial). The
adelic reformulation of the theory of automorphic forms provides a
natural spectral expansion of this space in which the basis of
Eisenstein series is indexed by a set of parameters of the
form 
\begin{equation}\label{eisenpara}
  \{(t, \chi_{1},\chi_{2}, b)\mid t \in \mathbb{R}, \chi_1\chi_{2}=\chi,\
b\in\mathcal{B}(\chi_{1},\chi_{2})\},
\end{equation}
where $\mathcal{B}(\chi_{1},\chi_{2})$ is some finite (possibly empty) set\footnote{specifically, $\mathcal{B}(\chi_{1},\chi_{2})$
corresponds to an orthonormal basis in the space of an induced
representation constructed out of the pair $(\chi_{1},\chi_{2})$, 
but we do not need to be more precise.}. We refer to \cite{GJ} for the
definition of these parameters as well as for the proof of the corresponding version of the
spectral theorem.   The reason why we introduce this particular basis of $\mathcal{E}(N, \chi)$  is that it consists of Hecke
eigenforms for all Hecke operators $T_n$ with $(n, N) = 1$ where the Hecke eigenvalue of the corresponding Eisenstein series  $E_{t, \chi_1, \chi_2, b}$ is given by 
$$\lambda_{t, \chi_{1},\chi_{2},b}(n)=\sum_{d_1d_2=n}\chi_{1}(d_1)d_1^{it}\chi_{2}(d_2)d_2^{-it}.$$
Of course, these Eisenstein series can be described in purely classical terms without recourse to adeles. If a cusp form $f$ is an eigenform of some Hecke operator $T_n$, we denote the corresponding eigenvalue by $\lambda_f(n)$. A form $f$ in $\mathcal{B}_k(N, \chi)$, $\mathcal{B}(N, \chi)$, or an Eisenstein series $f = E_{t, \chi_1, \chi_2, b}$ has a Fourier expansion
\begin{displaymath}
  f(z) = \rho_f(0, y) + \sum_{n \not= 0} \rho_f(n) W_f(|n| y)e(nx) 
\end{displaymath}
where
\begin{equation}\label{defW}
  W_f(y) = \left\{\begin{array}{ll}
     \Gamma(k)^{-1/2} (4 \pi y)^{k/2}e^{-2\pi y}, & f \in \mathcal{B}_k(N, \chi),\\
     y^{1/2} \cosh(\frac{1}{2}\pi t) K_{i t}(2 \pi y), & f \in \mathcal{B}(t, N, \chi) \text{ or } f = E_{t, \chi_1, \chi_2, b}
  \end{array}\right.
\end{equation}
for $y > 0$, 
\begin{equation}\label{sign}
  \rho(-n) = \eta \rho(n)
\end{equation}
for $n> 0$ and $\eta \in \{-1, 0, 1\}$, and $\rho_f(0, y)=0$ unless $f$ is Eisenstein.  The case $\eta = 0$ occurs for holomorphic forms. By our choice of (Hecke) bases for $\mathcal{S}_k(N, \chi), \mathcal{L}^2_0(N, \chi)$ and $\mathcal{E}(N, \chi)$ one has for any basis form $f$, any $n \geq 1$, and $(m,N)=1$  the relation
\begin{equation}\label{eig2}
 \lambda_{f}(m)\sqrt{n}\rho_f(n)  = \sum_{d\mid(m,n)}\chi(d)\sqrt{\frac{mn}{d^2}}\rho_{f}\left(\frac{mn}{d^2}\right);\\
\end{equation}
in particular, 
\begin{equation}\label{eig}\lambda_{f}(m)\rho_f(1)
=\sqrt{m}\rho_{f}({m})
\end{equation} 
for $(m,N)=1$. Moreover, these relations hold for all $m$, $n$ if $f$ is a cuspidal newform. For  such a cuspidal newform $f$ we have the uniform bounds \cite[Prop. 19.6]{DFI}
\begin{equation}\label{boundHecke}
  \sum_{n \leq X} |\lambda_f(n)|^2 \ll_{\varepsilon} X (XNt^{\ast})^{\varepsilon}
\end{equation}
for any $\varepsilon > 0, X\geq 1$, 
\begin{equation}\label{rho1}
 \frac{1}{N}(Nt^{\ast})^{-\varepsilon} \ll_{\varepsilon}  |\rho_f(1)|^2 \ll_{\varepsilon} \frac{1}{N} (Nt^{\ast})^{\varepsilon}
\end{equation}
(cf.\ \cite[(30)-(31)]{HM}) due to Iwaniec and Hoffstein-Lockhart, and
\begin{equation}\label{ind}
  |\lambda_f(n)| \ll_{\varepsilon} n^{\theta+\varepsilon}
\end{equation}
for $\theta = 7/64$, see \cite{KS}.  Note that we have normalized the Whittaker function $W_f$ in \eqref{defW} differently than in \cite[p.\ 593]{HM}. It follows from Atkin-Lehner theory \cite[Section 2]{AU} that
\begin{equation}\label{atkinlehner}
  \sup_{z \in \mathbb{H}} |f(z)| = \sup_{y \geq \frac{\sqrt{3}}{2N}} |f(x+iy)| 
\end{equation}
for a cuspidal newform $f$.  Here we use the fact that $N$ is square-free.\\ 

\section{Summation formulae}

First we recall the  Voronoi summation formulae   \cite[Appendix A3 and A4]{KMV}. Let $q\geq 1$ be an integer, and write $N_1 := (N, q)$, $N_2 = N/N_1$. Since $N$ is square-free, this implies $(N_1, N_2) = 1$.  For a cusp form $f$ and $y> 0$  define   $\mathcal{J}^{\pm} = \mathcal{J}_f^{\pm}$ by
\begin{equation}\label{defJ}
  \mathcal{J}^{+}(y) = 2\pi J_{k-1}(4\pi y), \quad \mathcal{J}^{-}(y) = 0
  \end{equation}
if $f$ is   holomorphic of weight $k$, and
\begin{equation}\label{defJ1}
  \mathcal{J}^{+}(y)  = \frac{\pi}{\cosh(\pi t)}\left(Y_{2i t}(4\pi y)+ Y_{-2 it}(4 \pi y)\right),\quad 
\mathcal{J}^{-}(y)  = 4\cosh(\pi t)K_{2i t}(4 \pi y)
 \end{equation}
if $f$ is Maa{\ss}. 

\begin{lemma} Let $(a, q) = 1$ and let $g$ be a smooth function, compactly supported in $(0, \infty)$. Let $f$ be a newform of level $N$. Then there exist two complex numbers $\eta_{\pm}$ of absolute value 1 (depending on $a, q$ and $f$) and a newform $f^{\ast}$ of the same level $N$ and the same archimedean parameter such that
\begin{equation}\label{voronoi}
  \sum_{n} \lambda_f(n) e\left(\frac{an}{q}\right) g(n) = \sum_{\pm} \frac{\eta_{\pm}}{q\sqrt{N_2}} \sum_{n} \lambda_{f^{\ast}}(n) e\left(\mp \frac{n\overline{aN_2}}{q}\right) \int_{0}^{\infty} g(\xi) \mathcal{J}^{\pm}\left(\frac{\sqrt{n\xi}}{q\sqrt{N_2}}\right)d\xi
\end{equation}
 where as usual the bar denotes the multiplicative inverse.
\end{lemma} 

For $m, n, c \in \Bbb{N}$ let $S_{\chi}(m, n; c)$ denote the twisted Kloosterman sum
\begin{displaymath}
  S_{\chi}(m, n; c) = \sum_{\substack{a \, (c)\\ (a, c) = 1}} \bar{\chi}(a) e\left(\frac{m \bar{a} + n a}{c}\right).
\end{displaymath}
Let $\phi : [0, \infty) \rightarrow \mathbb{C}$ be a smooth function
satisfying $\phi(0) = \phi'(0) = 0$, $\phi^{(j)}(y) \ll_\varepsilon
(1+y)^{-2-\varepsilon}$ for $0 \leq j \leq 3$. Let
\begin{equation}\label{besseltrans}
\begin{split}
  \dot{\phi}(k) &:= i^k\int_0^{\infty} J_{k-1}(x) \phi(y) \frac{dy}{y},\\
  \tilde{\phi}(t) & := \frac{i}{2\sinh(\pi t)}\int_0^{\infty} \left(J_{2it}(y) - J_{-2it}(y)\right)\phi(y)\frac{dy}{y}
\end{split}
\end{equation}
be Bessel transforms. 

\begin{lemma} For positive integers $m$, $n$ the trace
formula of Bruggeman--Kuznetsov holds:
\begin{displaymath}\begin{split} \sum_{c} \frac{1}{Nc}
S_{\chi}(m, n;Nc)\phi\left(\frac{4\pi\sqrt{mn}}{Nc}\right) =
&4\sum_{\substack{\kappa \geq 2 \text{ even}\\f\in \mathcal{B}_{\kappa}(N,\chi)}}
\dot{\phi}(k)
\overline{\rho_{f}(m)}\rho_{f}(n)\sqrt{mn}\\
+ \sum_{f\in\mathcal{B}(N,\chi)} 4\pi  \tilde{\phi}( t_f) 
\overline{\rho_{f}(m)}\rho_f(n)\sqrt{mn} \,
+ &\sum_{\substack{\chi_{1}\chi_{2}=\chi\\ b\in
\mathcal{B}(\chi_{1},\chi_{2})}}\int_{-\infty}^{\infty}\tilde{\phi}(t) \overline{\rho_{t, \chi_1, \chi_2, b}(m)}\rho_{t, \chi_1, \chi_2, b}(n)\sqrt{mn}\,dt.
\end{split}\end{displaymath}
\end{lemma}
For a proof see \cite[Theorem 9.8]{Iw1} and note the different normalization there\footnote{Note that equation \cite[(9.15)]{Iw1} should have the normalization factor $\frac{2}{\pi}$ instead of $\pi$, and in equation \cite[(B.49)]{Iw1} a factor 4 is missing.}  that is given by \cite[(1.26), (8.1), (8.2), (8.5), (8.6)]{Iw1}.  The proof of this formula with the 
Eisenstein parameters \eqref{eisenpara} is identical. The holomorphic counterpart   is
Petersson's trace formula 
\begin{lemma} We have
\begin{displaymath}
  \delta(m, n) + 2 \pi i^{-k} \sum_{c} \frac{1}{Nc} S_{\chi}(m, n; Nc)J_{k-1}\left(\frac{4\pi \sqrt{mn}}{Nc}\right)
  =\frac{4\pi\sqrt{mn}}{k-1} \sum_{f\in\mathcal{B}_{k}(N,\chi)}\overline{\rho_{f}(m)}\rho_{f}(n)
\end{displaymath}
where $\delta(m, n) = 1$ if $m=n$ and $\delta(m, n) = 0$ otherwise.
\end{lemma}
This is proved in \cite[Theorem 3.6]{Iw2} for Fourier coefficients normalized as in \cite[p.\ 52]{Iw2}.\\

Finally we need the following result. As usual, for $\alpha \in \Bbb{R}$ we write $\| \alpha \| := \min(|\alpha - n| : n \in \Bbb{Z})$. 
\begin{lemma} Let $0 < T \leq Z$, $\alpha \in \Bbb{R} \setminus \Bbb{Z}$. Let $\Phi$ be a smooth function, compactly supported on $[Z/2, 2Z]$ such that $\Phi^{(j)} \ll_j  T^{-j}$ for all $j \in \Bbb{N}_0$. Then
\begin{displaymath}
  \sum_{m} e(\alpha m) \Phi(m) \ll_j  Z (T\| \alpha\|)^{-j}
\end{displaymath}
for any $j \geq 2$. 
\end{lemma}

\textbf{Proof.} This is a simple application of Poisson summation and repeated integration by parts. 

 \section{Special functions}

For future reference we  record some properties of the functions $W_f$ defined in \eqref{defW} and $\mathcal{J}_f$ defined in \eqref{defJ} - \eqref{defJ1}.  For $y > 0$, 
$\sigma > 0$, $s = \sigma + it \in \Bbb{C}$, $k \in \mathbb{N}_0$, $A \geq 0$ and $\varepsilon > 0$ we have the  bounds
\begin{equation}\label{bessel3}
  J_{k}(y) \ll \frac{1+k}{1+y^{1/2}}, 
\end{equation}
\begin{equation}\label{bessel4}
  \frac{Y_{s}(y)}{\cosh(\pi t/2)} \ll_{\sigma, \varepsilon} \left\{\begin{array}{ll}
      ((1+|t|)y^{-1})^{|\sigma| + \varepsilon}, & y \leq 1 + |t|,\\
      ((1+|t|)y^{-1})^{- \varepsilon}, & 1+|t| < y \leq 1 + |s|^2,\\
      y^{-1/2}, & y > 1+ |s|^2,\\
   \end{array}\right.
  \end{equation}
 \begin{equation}\label{bessel5}
  \cosh(\pi t/2)K_{s}(y) \ll_{\varepsilon, \sigma} \left\{\begin{array}{ll} ((1+|t|)y^{-1})^{|\sigma|+\varepsilon}, & y \leq 1+ \frac{1}{2}\pi |t|\\ e^{-y + \pi|t|/2}y^{-1/2}, & y > 1 + \frac{1}{2}\pi|t| \end{array}\right\} \ll_A  \left(\frac{1+|t|}{y}\right)^{|\sigma| +\varepsilon} \left(1+\frac{y}{1+|t|}\right)^{-A}.
\end{equation}
The estimates \eqref{bessel3} - \eqref{bessel5} follow directly from 
\cite{HM}, Propositions 8 and 9. They are not best possible in all ranges, but sufficiently uniform for our present purposes. Moreover, we have
\begin{equation}\label{bessel1} 
\begin{split}
&  J_r(y)' =  \frac{1}{2}\left(J_{r-1}(y) - J_{r+1}(y)\right), \quad 
 K_r(y)' =  -\frac{1}{2}\left(K_{r-1}(y) + K_{r+1}(y)\right)
 \end{split}
\end{equation}
for arbitrary $r \in \mathbb{C}$. If  $\alpha > 0$, $g \in \mathcal{C}^1_0((0, \infty))$, then it follows from the formulas $(y^rJ_r(y))' = y^r J_{r-1}(y)$, $(y^rY_r(y))' = y^r Y_{r-1}(y)$, $(y^rK_r(y))' = -y^r K_{r-1}(y)$ and one integration by parts that 
\begin{equation}\label{bessel2}
  \int_{0}^{\infty} g(y)\left\{\begin{array}{lll}J_r(\alpha \sqrt{y})\\ Y_r(\alpha \sqrt{y})\\ K_r(\alpha \sqrt{y})\end{array} \right\} dy = \left\{\begin{array}{lll}+\\+\\-\end{array}\right\} \frac{2}{\alpha} \int_0^{\infty} \left(g'(y) \sqrt{y} - \frac{r}{2} \frac{g(y)}{\sqrt{y}}\right)\left\{\begin{array}{lll} J_{r+1}(\alpha\sqrt{y}) \\ Y_{r+1}(\alpha\sqrt{y})\\ K_{r+1}(\alpha\sqrt{y})\end{array}\right\} dy.
\end{equation}

 \begin{lemma} Let $f$ be a holomorphic cusp form or a non-exceptional Maa{\ss} cusp form (i.e.  $\Im t_f = 0$). Then one has
 \begin{displaymath}
  \frac{d^j}{dy^j} W_f(y) \ll_{j, \varepsilon, A} (t^{\ast})^{1/2} \left(\frac{t^{\ast}}{y}\right)^{j+\varepsilon}\left(1+\frac{y}{t^{\ast}}\right)^{-A}
\end{displaymath}
for   any $j, A \geq 0$.
\end{lemma}

\textbf{Proof.} If $f$ is Maa{\ss}, the lemma follows immediately from \eqref{defW}, the second equation in \eqref{bessel1},  and \eqref{bessel5}, and it conveniently simplifies further computations if we replace $y^{1/2}$ by $(t^{\ast})^{1/2}$. If $f$ is holomorphic, we observe that  the function
\begin{displaymath}
  W_f(y) = \Gamma(k)^{-1/2}  (4 \pi y)^{k/2} e^{-2 \pi y}
\end{displaymath}
 has its maximum at $y_0 = k/(4 \pi)$, and by Stirling's formula $W_f(y_0) \asymp k^{1/4}$. This readily confirms (a stronger version of) Lemma 5 for $j=0$ and $y \leq 10k$, say. For $y >10k$ we have again by Stirling's formula $W_f(y) \ll k^{1/4} \exp(-2\pi y + \frac{k}{2}\log\frac{4\pi e y}{k}) \leq k^{1/4} e^{-\pi y}$. 
 For $j > 0$, we observe that $W^{(j)}_f(y) \ll_j (1 + k/y)^j W_f(y)$ which completes the proof in the holomorphic case.  \\

\begin{lemma} Let $f$ be as in the preceding lemma. Let $1 \leq T \leq Z$, $\varepsilon > 0$ and let $g$ be a smooth function, compactly supported on $[Z/2, 2Z]$ such that $g^{(j)} \ll_j T^{-j}$ for all $j \in \Bbb{N}_0$. Let $\alpha > 0$ and 
\begin{displaymath}
  I := \int_0^{\infty} g(\xi) \mathcal{J}^{\pm}_f(\alpha \sqrt{\xi}) d\xi.
\end{displaymath}
Then
\begin{equation}\label{boundJ1}
  I \ll_{\varepsilon}  \frac{Z^{3/4} t^{\ast}}{\alpha^{1/2}}  (t^{\ast} (\alpha+\alpha^{-1}))^{\varepsilon}
\end{equation}
Moreover, if $\alpha \sqrt{Z/2} \geq 2 t^{\ast}$, then
\begin{equation}\label{boundJ2}
  I \ll _{\varepsilon, j} \left(\frac{1}{\alpha}\left(\frac{\sqrt{Z}}{T} + \frac{t^{\ast}}{\sqrt{Z}}\right)\right)^j  \frac{Z^{3/4} t^{\ast}}{\alpha^{1/2}}  (t^{\ast} (\alpha+\alpha^{-1}))^{\varepsilon}
\end{equation}
for any $j \in \Bbb{N}_0$.  
\end{lemma}

\textbf{Proof.} From \eqref{bessel3} - \eqref{bessel5} we infer $\mathcal{J}^{\pm}_f(\alpha \sqrt{\xi}) \ll_{\varepsilon} t^{\ast} (\alpha \sqrt{Z})^{-1/2} (t^{\ast}(\alpha+\alpha^{-1}))^{\varepsilon}$ for $\xi \in [Z/2, 2Z]$, and the first bound is immediate. For the proof of the second bound we apply $j$ times the formula \eqref{bessel2}  and then estimate trivially using \eqref{bessel3} - \eqref{bessel5} as before. \\

Finally we specify a test function $\phi$ for the Kuznetsov formula (Lemma 2).   For integers $2 \leq B < A$ of the same parity let us define
\begin{equation}\label{defphi}
  \phi_{A, B}(x) := i^{B-A} J_A(x)x^{-B}.
\end{equation}
Using
\cite[6.574.2]{GR} it is straightforward to verify that
\begin{equation}\label{phitrafo}
\begin{split}
  &\dot{\phi}_{A, B}(k) = \frac{B!}{2^{B+1}\pi}
  \prod_{j=0}^B \left\{\left(\frac{(1-k)i}{2}\right)^2 + \left(\frac{A+B}{2}-j\right)^2\right\}^{-1}
  \asymp_{A, B}\, \pm \,k^{-2B-2},\\
  &\tilde{\phi}_{A, B}(t) = \frac{B!}{2^{B+1} \pi}
  \prod_{j=0}^{B} \left\{t^2 + \left(\frac{A+B}{2}-j\right)^2\right\}^{-1}
  \asymp_{A, B}\, (1+|t|)^{- 2B-2}.
\end{split}
\end{equation}
In
particular,
\begin{equation}
\begin{split}\label{posdot}
\dot{\phi}_{A, B}(k) > 0 \quad &\text{for}\quad 2\leq k \leq A-B,\\
\tilde{\phi}_{A, B}(t)> 0\quad &\text{for}\quad t \in \mathbb{R} \cup i[-\theta, \theta]. 
\end{split}
\end{equation}

\section{The major arcs}

We keep the notation developed so far. For the proof of Theorem 1 we can assume that $N$ is sufficiently large. It will be convenient to introduce  a smooth partition of unity: Let $Z= 2^{\nu}$, $\nu \in \mathbb{N}_0$,  run through powers of 2, and let $G$ be a  smooth nonnegative function $G$ supported on $[1/2, 2]$ such that $\sum_Z G(x/Z) = 1$ for all $x \geq 1$.  For a (non-exceptional) newform $f \in \mathcal{B}(N, \chi)$ or $\mathcal{B}_k(N, \chi)$ let 
\begin{displaymath}
  f^{Z}(x+iy) := \sum_{n} \rho_f(n) W_f(ny) G\left(\frac{n}{Z}\right) e(nx). 
\end{displaymath}
 In view of \eqref{atkinlehner} we can assume $y \geq 1/(2N)$. For notational simplicity let
\begin{displaymath}
  P := t^{\ast}N.
\end{displaymath}
First we record the trivial bound
\begin{equation}\label{trivial1}
    |f^Z(x+iy)| \ll_{A, \varepsilon}   \frac{P^{\varepsilon} (t^{\ast}Z)^{1/2}}{N^{1/2}} \left(1+ \frac{Zy}{t^{\ast}}\right)^{-A}
\end{equation}
for any $A \geq 0$ and any $\varepsilon> 0$. This follows directly from Cauchy-Schwarz in combination with \eqref{eig} - \eqref{rho1} and Lemma 5 with $j=0$. Fix $\varepsilon > 0$ and let 
\begin{equation}\label{boundZ}
  Z_0 := P^{\varepsilon}t^{\ast}/y \leq 2P^{\varepsilon} t^{\ast} N.
\end{equation}
Then we have
\begin{lemma} If $\varepsilon > 0$, then
\begin{displaymath}
  |f(x+iy)| \ll_{\varepsilon} P^{\varepsilon} \max_{Z \leq Z_0} |f^{Z}(x+iy)|,
  \end{displaymath}
and
\begin{equation}\label{trivialZ}
    |f^Z(x+iy)| \ll_{\varepsilon}   \frac{P^{\varepsilon} (t^{\ast}Z)^{1/2}}{N^{1/2}}. \end{equation}
In particular, 
\begin{equation}\label{trivial}
  |f(x+iy)| \ll_{\varepsilon} t^{\ast}(Nt^{\ast})^{\varepsilon}.
\end{equation}      
\end{lemma}
Indeed, choosing $A$ in \eqref{trivial1} sufficiently large in terms of $\varepsilon$, it is obvious that $\sum_{Z = 2^{\nu} \geq Z_0} |f^Z(x + iy)|$ is negligible, and by \eqref{sign} it is enough to look at positive Fourier coefficients. This establishes the first bound, while \eqref{trivialZ} and \eqref{trivial} are direct consequences of \eqref{trivial1} and \eqref{boundZ}. 

We remark that Lemma 7 implies already a result of the type of Theorem 1 as soon as $Z_0 \leq N^{1-\delta}$, so the critical range of the dyadic decomposition  is really $Z= N^{1+o(1)}$. In Proposition 1 below we will prove a  version of \eqref{trivial} that is stronger in the $t^{\ast}$-aspect. We observe that we can assume 
\begin{equation}\label{ranget}
  t^{\ast} \leq N^{1/165},
\end{equation}
for otherwise Theorem 1 follows from \eqref{trivial}. 


The aim of this section is to improve on the trivial bound  if $x$ can be ``well-approximated". This is summarized in the following lemma.
\begin{lemma} Let  $f$ be a newform as above and write
\begin{displaymath}
  x = \frac{a}{q} + \beta
\end{displaymath}
where $(a, q) = 1$, $q \leq N$  and $\beta \leq 1$, say.  Let $N_2 := N/(N, q)$ as in Lemma 1. Then
\begin{displaymath}
  |f^Z(x+iy)| \ll_{\varepsilon} P^{\varepsilon}(t^{\ast})^{3/2}q\left(\beta^{3/2}Z + \frac{(t^{\ast})^{3/2}}{Z^{1/2}}\right).
\end{displaymath} 
\end{lemma}

\textbf{Proof.} Using \eqref{rho1} and the Voronoi summation formulae \eqref{voronoi} we find by trivial estimation that
\begin{equation}\label{wefind}
\begin{split}
|f^Z(x+iy)|& = |\rho_f(1)| \left| \sum_{n} \lambda_f(n) e\left(n\frac{a}{q}\right)\frac{G(n/Z)W_f(ny)e(n\beta)}{\sqrt{n}}\right|\\
&  \ll_{\varepsilon} \frac{(Nt^{\ast})^{\varepsilon}}{\sqrt{N}q\sqrt{N_2}}  \sum_{\pm} \sum_{n}  |\lambda_{f^{\ast}}(n)|\left|\int_0^{\infty} \frac{G(\xi/Z)W_f(\xi y)e(\xi\beta)}{\sqrt{\xi}} \mathcal{J}^{\pm}\left(\frac{\sqrt{n\xi}}{q\sqrt{N_2}}\right)d\xi\right|. 
\end{split}
\end{equation}
Writing $g(\xi) := G(\xi/Z)W_f(y\xi)e(\xi\beta)\xi^{-1/2}$, Lemma 5 implies
\begin{displaymath} 
   g^{(j)}(\xi) \ll_{j, \varepsilon} P^{\varepsilon} \left(\frac{t^{\ast}}{ Z}\right)^{1/2} \left(\frac{t^{\ast}}{Z} + \beta \right)^{j}
\end{displaymath}   
    for any $j \in \Bbb{N}_0$, $\varepsilon > 0$.  Let us first assume
\begin{equation}\label{X}
  n \geq X :=  8 P^{\varepsilon} q^2 N_2 \left(\beta^2Z+\frac{(t^{\ast})^2}{Z}\right).
\end{equation}
Then $\sqrt{n\xi}/q\sqrt{N_2} \geq 2P^{\varepsilon} t^{\ast}$, and hence by \eqref{boundJ2} the integral in the second line of  \eqref{wefind} is at most
\begin{equation}\label{typicaltrick}
  \ll_{\varepsilon, j} P^{\varepsilon} \left(\frac{t^{\ast}}{ Z}\right)^{1/2}  \left(\frac{q\sqrt{N_2}}{\sqrt{n}}\left(\beta\sqrt{Z}+ \frac{t^{\ast}}{\sqrt{Z}}\right)\right)^{j} Z^{3/4} t^{\ast}\left(\frac{\sqrt{n}}{q\sqrt{N_2}}\right)^{-1/2} \ll P^{20} P^{-(j-2) \varepsilon} n^{-3/2}
\end{equation}
for $1/2 \leq Z \leq Z_0$. Choosing $j$ sufficiently large (in terms of $\varepsilon$), we see that these values of $n$ contribute a negligible amount to \eqref{wefind}. For the remaining $n$ we use \eqref{boundJ1}, and bound \eqref{wefind} by 
\begin{displaymath}
  \ll_{\varepsilon} \frac{P^{\varepsilon}}{\sqrt{N}q\sqrt{N_2}}   \sum_{n \leq X}  |\lambda_{f^{\ast}}(n)|  \left(\frac{t^{\ast}}{ Z}\right)^{1/2}  Z^{3/4} t^{\ast}\left(\frac{\sqrt{n}}{q\sqrt{N_2}}\right)^{-1/2}. 
\end{displaymath}
We perform the sum over $n$ by Cauchy-Schwarz and \eqref{boundHecke} obtaining
\begin{displaymath}
  \ll_{\varepsilon} \frac{P^{\varepsilon} (t^{\ast})^{3/2} Z^{1/4}X^{3/4}}{N^{1/2} q^{1/2} N_2^{1/4}}.
\end{displaymath}
Now we use \eqref{X} together with $N_2 \leq N$ and Lemma 8 is immediate. \\

\emph{Remark.} Obviously this estimate will only be useful if $\beta$ and $q$ are relatively small, i.e. $x$ is well-approximable. For example, if $Z \approx N$, Lemma 8 is non-trivial under assumptions a little stronger than  $q \leq N^{1/2}$ and $q \beta^{3/2} \leq 1/N$ which defines a subset of $[0, 1]$ of measure $N^{-\delta}$.

\section{Amplification}

In this section we prepare the scene for a different method to estimate  $|f^Z(z)|=|f^Z(x+iy)|$ that turns out to be strong if $x$ is off rational points with small denominators. Essentially, we will consider an amplified second moment. More precisely, let $f_0$ be the newform whose $L^{\infty}$-norm we want to bound. Let 
\begin{equation}\label{rangeL}
  \max(N^{1/100}, (t^{\ast})^3) \leq L \leq N
\end{equation}  
be a parameter   to be chosen later, let 
\begin{displaymath}
  \Lambda_1 := \{p \nmid N : p \in [L, 2L] \text{ prime}\}, \quad \Lambda_2 := \{p^2 : p \in \Lambda_1\}
\end{displaymath}
and let us define
\begin{equation}\label{ampli}
  \alpha(\ell) := \left\{\begin{array}{ll}
  \lambda_{f_0}(\ell)\bar{\chi}(\ell), &  \ell \in \Lambda_1,\\
  -\bar{\chi}(\ell), &   \ell \in \Lambda_2,\\
  0, & \text{else.}
  \end{array}
   \right.
\end{equation}
For $\kappa \in \{2, 4, 6, \ldots\}$ let us define 
\begin{displaymath}
  \mathcal{Q}_{\kappa}(z) := \frac{2}{\kappa-1} \sum_{f \in \mathcal{B}_{\kappa}(N, \chi)} \left| \sum_{\ell} \lambda_f(\ell) \alpha(\ell)\right|^2 |f^Z(z)|^2
\end{displaymath}
and 
\begin{displaymath}
\begin{split}
  \mathcal{Q}(z) := &\sum_{\kappa \geq 2 \text{ even}} 2(\kappa - 1)\dot{\phi}(\kappa) \mathcal{Q}_{\kappa}(z) \\
    & + 4\pi \sum_{f \in \mathcal{B}(N, \chi)} \tilde{\phi}(t_f)   \left| \sum_{\ell} \lambda_f(\ell) \alpha(\ell)\right|^2 \left| \sum_{n} \rho_f(n) G(n/Z)W_{f_0}(ny)e(nx)\right|^2\\
    & + \sum_{\substack{\chi_1 \chi_2 = \chi\\ b \in \mathcal{B}(\chi_1, \chi_2)}} \int_{\mathbb{R}} \tilde{\phi}(t)  \left| \sum_{\ell} \lambda_{t, \chi_1, \chi_2, b}(\ell) \alpha(\ell)\right|^2 \left| \sum_{n} \rho_{t, \chi_1, \chi_2, b}(n) G(n/Z) W_{f_0}(ny)e(nx)\right|^2 dt.
    \end{split}
\end{displaymath}
Here $\phi$ is the function defined in \eqref{defphi} with some small $B\leq 5$ and some very large $A$ to be chosen later, and $\dot{\phi}$ and $\tilde{\phi}$ are the integral transforms \eqref{besseltrans} of $\phi$ that satisfy the relations \eqref{phitrafo} and \eqref{posdot}.  Using  Lemma 7 and \eqref{phitrafo}, we see that the $\kappa$-sum in $\mathcal{Q}(z)$ is absolutely convergent.   The quantity $\mathcal{Q}(z)$ is not really an amplified second moment, but rather a kind of ``fake moment", since we consider truncated cusp forms and we need the weight function $W_{f_0}$ to be independent of $f$. Finally we remark that the amplifier that we chose is most likely not the most efficient, but it seems hard to show that the more natural amplifier $\alpha(\ell) = \lambda_{f_0}(\ell)$ for all $\ell \in [L, 2L]$ gives a nontrivial result. If we knew $\sum_{\ell \sim L} |\lambda_{f_0}(\ell)|^2 \gg L^{1-\varepsilon}$, then our analysis would no longer involve those $\ell \sim L^2$. This would greatly simplify the following sections and improve the exponent in Theorem 1 by about a factor 2.

 From \eqref{eig2} and \eqref{eig} we conclude $\lambda_f(p)^2 = \lambda_f(p^2) + \chi(p)$, hence
\begin{equation}\label{sizeampli}
  \sum_{\ell} \lambda_{f_0}(\ell)\alpha(\ell) = \sum_{p \in \Lambda_1} 1 \gg L N^{-\varepsilon}.
\end{equation}
Dropping all but one term in the definitions of $\mathcal{Q}_\kappa$ and $\mathcal{Q}(z)$ and using \eqref{phitrafo}, \eqref{posdot} and \eqref{sizeampli}, we see that 
\begin{equation}\label{seethat}
  |f^Z_0(z)|^2 \ll_{\varepsilon}   P^{\varepsilon} L^{-2} t^{\ast} \mathcal{Q}_k(z)
\end{equation}
if $f_0$ is a newform in $\mathcal{S}_k(N, \chi)$ for even $k > A-5$, and 
\begin{equation}\label{seethat2}
  |f^Z_0(z)|^2 \ll_{\varepsilon}  P^{\varepsilon} L^{-2} (t^{\ast})^{2B+2} \left(\mathcal{Q}(z) + \sum_{\kappa \geq A-4} \kappa|\dot{\phi}(\kappa)| \mathcal{Q}_{\kappa}(z)\right)
\end{equation}
if $f_0$ is a newform in $\mathcal{S}_k(N, \chi)$ for even $k \leq A-5$, or  a newform in $\mathcal{L}^2_0(N, \chi)$. We proceed to bound $\mathcal{Q}_{\kappa}(z)$ and $\mathcal{Q}(z)$ from above. In view of \eqref{seethat} and \eqref{seethat2}, we need to get a bound of the form $L^{2-\delta}$ for some $\delta > 0$ and some sufficiently small $L$ in terms of the level $N$. Opening the squares, using\footnote{Note that this relation holds in particular for Eisenstein series}  \eqref{eig2}, and applying the appropriate trace formula Lemma 2 or Lemma 3, we obtain
\begin{equation}\label{weobtain1}
\begin{split}
 \mathcal{Q}_{\kappa}(z)&  =  \sum_{\ell_1, \ell_2} \bar{\alpha}(\ell_1)\alpha(\ell_2) \sum_{n, m} \frac{G(m/Z)\bar{W}_{f_0}(my) G(n/Z)W_{f_0}(ny)e((n-m)x)}{\sqrt{nm}} \sum_{d_1 \mid (m, \ell_1)}\sum_{d_2 \mid (n, \ell_2)} \bar{\chi}(d_1)\chi(d_2)\\
  &\quad\quad\quad\quad \times  \left(\frac{1}{2\pi}\delta\left(\frac{m\ell_1}{d_1^2}, \frac{n\ell_2}{d_2^2}\right) +  \sum_c \frac{i^{-\kappa}}{Nc} S_{\chi}\left(\frac{m\ell_1}{d_1^2}, \frac{n\ell_2}{d_2^2}; Nc\right) \phi\left(\frac{4\pi \sqrt{nm\ell_1\ell_2}}{d_1d_2Nc}\right)\right)\\
    & =:  \mathcal{Q}^{\Delta}_{\kappa}(z)  + \mathcal{Q}^{\text{off-}\Delta}_{\kappa}(z),
 \end{split} 
\end{equation}
say, where $\phi = J_{\kappa -1} = i^{\kappa-1}\phi_{\kappa-1, 0}$ in the notation \eqref{defphi}; and
\begin{displaymath}
\begin{split}
 \mathcal{Q}(z) & =  \sum_{\ell_1, \ell_2} \bar{\alpha}(\ell_1)\alpha(\ell_2) \sum_{n, m} \frac{G(m/Z)\bar{W}_{f_0}(my) G(n/Z)W_{f_0}(ny)e((n-m)x)}{\sqrt{nm}} \sum_{d_1 \mid (m, \ell_1)}\sum_{d_2 \mid (n, \ell_2)} \bar{\chi}(d_1)\chi(d_2)\\
  &\quad\quad\quad\quad\times     \sum_c \frac{1}{Nc} S_{\chi}\left(\frac{m\ell_1}{d_1^2}, \frac{n\ell_2}{d_2^2}; Nc\right) \phi\left(\frac{4\pi \sqrt{nm\ell_1\ell_2}}{d_1d_2Nc}\right) 
 \end{split} 
\end{displaymath}
where  $\phi =\phi_{A, B}$.  We are assuming that $\kappa$ and $A$ are sufficiently large, so that in either case $\phi$ decays sufficiently rapidly at $0$. The diagonal term $\mathcal{Q}^{\Delta}_{\kappa}(z) $ in \eqref{weobtain1} can be estimated trivially, and we summarize the result in the following lemma.

\begin{lemma} We have
\begin{displaymath}
 \mathcal{Q}^{\Delta}_{\kappa}(z)  \ll_{\varepsilon} P^{\varepsilon} t^{\ast} L
\end{displaymath}
for any $\varepsilon > 0$. 
\end{lemma}
To see this, we estimate
\begin{displaymath}
\begin{split}
 \mathcal{Q}^{\Delta}_{\kappa}(z)  &\leq \sum_{d_1, d_2} \sum_{\substack{\ell_1, \ell_2, n, m\\ m\ell_1 = n\ell_2}} \frac{|\alpha(d_1\ell_1)\alpha(d_2\ell_2) G(md_1/Z)W_{f_0}(md_1y)G(nd_2/Z)W_{f_0}(nd_2y)|}{\sqrt{d_1d_2nm}}\\
 & \ll \sum_{d_1, d_2, \ell_1, \ell_2} \frac{|\alpha(d_1\ell_1)\alpha(d_2\ell_2)| (\ell_1, \ell_2)}{\sqrt{d_1d_2\ell_1\ell_2}} \sum_m \frac{1}{m} \left|W_{f_0}\left(m\frac{\ell_2d_1}{(\ell_1, \ell_2)} y\right)W_{f_0}\left(m\frac{\ell_1d_2}{(\ell_1, \ell_2)}y\right)\right|\\
 & \ll_{\varepsilon} \sum_{d_1, d_2, \ell_1, \ell_2} \frac{|\alpha(d_1\ell_1)\alpha(d_2\ell_2)| (\ell_1, \ell_2)}{\sqrt{d_1d_2\ell_1\ell_2}} t^{\ast} P^{\varepsilon} 
\end{split}
\end{displaymath}
using Lemma 5. By \eqref{boundHecke} and the fact that $\alpha(\ell)$ is supported on prime powers (see \eqref{ampli}), we arrive at the bound of the lemma.

\section{Bounding exponential sums}

It remains to estimate the two similar quantities $\mathcal{Q}_{\kappa}^{\text{off-}\Delta}(z)$ and $\mathcal{Q}(z)$, an upper bound for which is given by  
\begin{displaymath}
\begin{split}
  \sum_{d_1, d_2, \ell_1, \ell_2} &\frac{|\alpha(d_1\ell_1)\alpha(d_2\ell_2)|}{\sqrt{d_1d_2}}\sum_c \frac{1}{Nc}\left| \sum_{n, m} \ S_{\chi}(m\ell_1, n\ell_2; Nc) \phi\left(\frac{4\pi\sqrt{nm\ell_1\ell_2}}{Nc}\right)\right.\\
  &\times\left.  \frac{G(d_1m/Z)\bar{W}_{f_0}(d_1my) G(d_2n/Z)W_{f_0}(d_2ny)e((d_2n-d_1m)x)}{\sqrt{nm}}  \right|.
\end{split}
\end{displaymath}
This is the most technical part of the paper and will be completed once we reach Lemma 11 below. Hecke eigenvalues have now disappeared from the picture (except from the amplifier), and we are left with bounding weighted exponential sums. We write
\begin{displaymath}
  x = \frac{u}{N} + \beta
\end{displaymath}
where $u \in \mathbb{Z}$ (not necessarily coprime to $N$) and $|\beta| \leq 1/N$.  For notational simplicity let us temporarily write 
\begin{equation}\label{defX}
  X_1 := Z/d_1, \quad X_2 := Z/d_2
\end{equation}
which corresponds to the sizes of $m$ and $n$ respectively. By the rapid decay of $\phi$ near 0, we can assume that 
\begin{equation}\label{ranges}
 c \leq  C:= \frac{\sqrt{X_1X_2\ell_1\ell_2}}{N}P^{\varepsilon}
  \end{equation}
for any $\varepsilon> 0$, for otherwise we can save as many powers of $P$ as we wish; this is a similar argument as in \eqref{typicaltrick}.  Opening the Kloosterman sum, we are left with
\begin{displaymath}
\begin{split}
     \sum_{\substack{d_1, d_2\\ \ell_1, \ell_2}}&\frac{|\alpha(d_1\ell_1)\alpha(d_2\ell_2)|}{\sqrt{d_1d_2}}\sum_{c \leq C}  \sum_{\substack{a \, (Nc)\\(a, Nc) = 1}} \frac{1}{Nc} \left|\sum_{n, m} e\left(m \left(\frac{\ell_1 \bar{a}}{Nc} - \frac{d_1u}{N}\right) + n \left(\frac{\ell_2a}{Nc} + \frac{d_2u}{N}\right)\right) \Phi_{\substack{d_1, d_2\\ \ell_1, \ell_2}}(m, n; Nc; y, \beta) \right|
\end{split}
\end{displaymath}
where 
\begin{displaymath}
  \Phi_{\substack{d_1, d_2\\ \ell_1, \ell_2}}(m, n; Nc; y, \beta) := \frac{\bar{W}_{f_0}(d_1my) W_{f_0}(d_2ny)G(m/X_1)G(n/X_2)e((d_2n-d_1m)\beta)}{\sqrt{nm}}  \phi\left(\frac{4\pi\sqrt{nm\ell_1\ell_2}}{Nc}\right).
\end{displaymath}
Recall that $\phi = \phi_{A, B}$ as in \eqref{defphi} for some fixed $A, B$  or $\phi = J_{\kappa -1}$ (where $A$ and $\kappa$ are very large). The general idea is now to observe that $\Phi$ oscillates only mildly, and so the $n, m$-sum can only be large if both $a$ and $\bar{a}$ are in a certain small range. This will happen rarely and, as it turns out, only if $u/N$ can be well-approximated in which case we can use the bound from Section 5. More precisely, in the range \eqref{ranges} and $Z \leq Z_0$ as in \eqref{boundZ} we have  by Lemma 5, \eqref{bessel1} and  \eqref{bessel3}
\begin{equation}\label{boundcapPhi}
 X_1^iX_2^j \frac{\partial^i}{\partial \xi^i} \frac{\partial^j}{\partial \eta^j}  \Phi_{\substack{d_1, d_2\\ \ell_1, \ell_2}}(\xi, \eta; Nc; y, \beta) \ll_{i, j, \varepsilon} \frac{t^{\ast} (Nc)^{1/2} P^{\varepsilon}}{(X_1X_2)^{3/4} (\ell_1\ell_2)^{1/4}} \left(t^{\ast} + \frac{\sqrt{X_1X_2\ell_1\ell_2}}{Nc}\right)^{i+j}
\end{equation}
for any $\varepsilon > 0$ and $i, j \in \mathbb{N}_0$ (recall $|\beta| \leq 1/N)$. The implied constant depends also on the choice of $A$ and $\kappa$ in $\phi$, but we regard these as fixed. Applying Lemma 4 to the $n$- and $m$-sum we see that the double sum is negligible, unless both
\begin{equation}\label{unlessboth}
\begin{split}
  \left\| \frac{\ell_1 \bar{a}}{Nc} - \frac{d_1u}{N} \right\| \leq  \frac{P^{\varepsilon}(t^{\ast} + \sqrt{X_1X_2\ell_1\ell_2}/(Nc))}{X_1}, \quad \left\|  \frac{\ell_2a}{Nc} + \frac{d_2u}{N} \right\| \leq   \frac{P^{\varepsilon}(t^{\ast} + \sqrt{X_1X_2\ell_1\ell_2}/(Nc))}{X_2},
\end{split}
\end{equation}
for otherwise we can save as many powers of $P$ as we wish by choosing $j$ sufficiently large in Lemma 4. If \eqref{unlessboth} is satisfied, we cannot do better than estimating the $n, m$-sum trivially by
\begin{displaymath}
  \ll \frac{t^{\ast} (Nc)^{1/2} (X_1X_2)^{1/4}P^{\varepsilon}}{ (\ell_1\ell_2)^{1/4}},
\end{displaymath}
using \eqref{boundcapPhi} with $i=j=0$. Now \eqref{unlessboth} implies 
\begin{equation}\label{congruence}
 \ell_1 \bar{a} \equiv d_1uc + r_1\, (Nc), \quad -\ell_2 a \equiv d_2 uc + r_2 \, (Nc)  
\end{equation}
for certain integers 
\begin{equation}\label{ranges2}
\begin{split}
 &  |r_1| \leq R_1 := P^{\varepsilon}\left(t^{\ast} + \frac{\sqrt{X_1X_2\ell_1\ell_2}}{Nc}\right) \frac{Nc}{X_1} \leq 2P^{2\varepsilon} t^{\ast}\left(\frac{\ell_1\ell_2d_1}{d_2}\right)^{1/2},\\
 &   |r_2| \leq R_2 := P^{\varepsilon}\left(t^{\ast} + \frac{\sqrt{X_1X_2\ell_1\ell_2}}{Nc}\right) \frac{Nc}{X_2} \leq 2P^{2\varepsilon} t^{\ast}\left(\frac{\ell_1\ell_2d_2}{d_1}\right)^{1/2},
  \end{split}
\end{equation}
cf.\ \eqref{defX} and \eqref{ranges}. 
If $d_1, d_2, \ell_1, \ell_2, r_1, r_2, c$ are given, $a$ is determined modulo $Nc/(c, \ell_1, \ell_2)$, hence the number of solutions $a$ modulo $Nc$ to \eqref{congruence} is at most $(c, \ell_1, \ell_2)$. Moreover, \eqref{congruence} implies  
\begin{equation}\label{cong}
  (d_1uc + r_1)(d_2uc+r_2) +\ell_1\ell_2 \equiv 0 \, (Nc)
\end{equation}
as well as
\begin{equation}\label{gcd}
  (\ell_1, c) \mid r_1, \quad (\ell_2, c) \mid r_2.
\end{equation}
Summarizing the preceding discussion, we obtain the following preliminary bound for the off-diagonal term 
\begin{equation}\label{atmost}
P^{\varepsilon} \sum_{\substack{d_1, d_2\\ \ell_1, \ell_2}} \frac{|\alpha(d_1\ell_1)\alpha(d_2\ell_2)|}{\sqrt{d_1d_2}}  \sum_{\substack{c \leq C, \, |r_1| \leq R_1, \, |r_2| \leq R_2\\ Nc \mid (d_1uc+r_1)(d_2uc+r_2) +\ell_1\ell_2\\ (\ell_1, c) \mid r_1,\, (\ell_2, c) \mid r_2}} \frac{(c, \ell_1, \ell_2)t^{\ast} Z^{1/2} }{(Nc)^{1/2} (d_1 d_2 \ell_1\ell_2)^{1/4}}.
 \end{equation} 
We will now split this expression depending on the sizes of $d_1, d_2, \ell_1, \ell_2 \in \Lambda_0 \cup \Lambda_1 \cup \Lambda_2$ where $\Lambda_0 := \{1\}$, cf.\ \eqref{ampli}.  For notational simplicity let us define
\begin{equation}\label{sizes}
C^{(j)} := \frac{P^{\varepsilon}Z (2L)^j}{N\sqrt{d_1 d_2}}, \quad R_1^{(j)} :=  2P^{\varepsilon} t^{\ast}(2L)^j \left(\frac{d_1}{d_2}\right)^{1/2}, \quad R_2^{(j)} :=  2P^{\varepsilon} t^{\ast}(2L)^j \left(\frac{d_2}{d_1}\right)^{1/2}. 
\end{equation}
Bounding the amplifier \eqref{ampli} by \eqref{ind} and cutting the $c$-sum into dyadic blocks, we can estimate  \eqref{atmost} as
\begin{equation}\label{gettingcloser}
\begin{split}
  \frac{P^{\varepsilon}t^{\ast} Z^{1/2} }{N^{1/2}}&  \sum_{\substack{i_1, i_2, j_1, j_2 \in \{0, 1, 2\}\\ 1 \leq i_1 + j_1 \leq 2\\ 1 \leq i_2 + j_2 \leq 2}}
   \sum_{\substack{d_1 \in \Lambda_{i_1}\\ d_2 \in \Lambda_{i_2}}} \frac{L^{\theta(4 - i_1- j_1 - i_2 -j_2)} }{(d_1d_2)^{3/4}L^{j/2}}  \max_{C \leq C^{(j)}} \frac{1}{C^{1/2}}  \sum_{\substack{C \leq c < 2C, \, \ell_1 \in \Lambda_{j_1}, \, \ell_2 \in \Lambda_{j_2} \\ |r_1| \leq R_1^{(j)}, \, |r_2| \leq R_2^{(j)}\\   Nc \mid (d_1uc+r_1)(d_2uc+r_2) +\ell_1\ell_2\\ (\ell_1, c) \mid r_1,\, (\ell_2, c) \mid r_2}}(c, \ell_1, \ell_2)
  \end{split}
\end{equation}
where $j := (j_1+j_2)/2$.  Note that $C$ is now independent of $\ell_1$ and $\ell_2$. We proceed to analyze the innermost (multiple) sum.  Considering \eqref{cong} modulo $c$, we find 
\begin{equation}\label{gleichung}
 r_1r_2 + \ell_1 \ell_2  = sc
\end{equation}  
   for some integer $s$ of size
\begin{equation}\label{sizeS}
  |s| \leq  \frac{R_1^{(j_1)}R_2^{(j_2)} + (2L)^{j_1+j_2}}{C} \leq \frac{5P^{2\varepsilon} (t^{\ast})^2 (2L)^{2j}}{C} =: S^{(j)} 
\end{equation}
by \eqref{ranges2} and \eqref{sizes}. Substituting this back into \eqref{cong}, we obtain
\begin{displaymath}
  u^2d_1d_2c + u(d_1r_2 + d_2r_1) +s \equiv 0 \, (N). 
\end{displaymath}
If $r_1, r_2, s, c$ are given, $\ell_1\ell_2$ is determined, and hence (up to allowed permutations of the prime factors) $\ell_1$ and $\ell_2$ are fixed. If $v_p$ denotes the usual $p$-adic valuation, then  \eqref{gleichung} implies for every prime $p$ the inequality  $v_p(s) \geq \min(v_p(\ell_1) + v_p(\ell_2)- v_p(c), v_p(r_1)+ v_p(r_2) - v_p(c))$ which together with \eqref{gcd} gives 
\begin{equation}\label{vps}
  v_p(s) \geq \min(v_p(\ell_1) + v_p(\ell_2) - v_p(c), v_p(\ell_1), v_p(\ell_2), v_p(c)). 
  \end{equation}
Since  $c \leq 8P^{2\varepsilon} t^{\ast} (\ell_1\ell_2)^{1/2}$ by \eqref{sizes} and \eqref{boundZ} and since $\ell_1, \ell_2$ are supported on prime powers, we obtain  by \eqref{rangeL} that $v_p(c) \leq \frac{1}{2}(v_p(\ell_1) + v_p(\ell_2))$ if $\varepsilon$ was chosen sufficiently small.  Therefore \eqref{vps} and \eqref{gcd} imply that
\begin{displaymath}
 (c, \ell_1, \ell_2) \leq (c, s, r_1, r_2). 
\end{displaymath}
We will see later that the case $(j_1, j_2) = (2, 2)$ in \eqref{gettingcloser} requires extra work. If $(j_1, j_2) \not= (2, 2)$ we replace the innermost sum of \eqref{gettingcloser} by
\begin{equation}\label{case1}
   \sum_{\substack{c \asymp C, \, |s| \leq S^{(j)} \\ |r_1| \leq R_1^{(j)}, \, |r_2| \leq R_2^{(j)}\\   N \mid u^2d_1d_2 c + u(d_1r_2 + d_2r_1) + s}}(c, s, r_1, r_2) \leq   \sum_{f \leq \min(C, S^{(j)}, R_1^{(j)}, R_2^{(j)})} f \sum_{\substack{c \asymp Cf^{-1}, \, |s| \leq S^{(j)}f^{-1} \\ |r_1| \leq R_1^{(j)}f^{-1}, \, |r_2| \leq R_2^{(j)}f^{-1}\\   N \mid u^2d_1d_2 c + u(d_1r_2 + d_2r_1) + s}}1.
\end{equation}
In the last step we used that $(\ell_1\ell_2, N) =1$ and \eqref{gleichung} imply  that $(c, s, r_1, r_2)$ is coprime to $N$.  If $(j_1, j_2) = (2, 2)$, we will not be able to bound \eqref{case1} satisfactorily, but we recall that in this case $\ell_1\ell_2$ is a perfect square, hence we only need to bound the quantity
\begin{equation}\label{case2}
  \sum_{\substack{c \asymp C, \, |s| \leq S^{(j)} \\ |r_1| \leq R_1^{(j)}, \, |r_2| \leq R_2^{(j)}\\  sc-r_1r_2  = \square\\  N \mid u^2d_1d_2 c + u(d_1r_2 + d_2r_1) + s}}(c, s, r_1, r_2) \leq   \sum_{f \leq \min(C, S^{(j)}, R_1^{(j)}, R_2^{(j)})} f \sum_{\substack{c \asymp Cf^{-1}, \, |s| \leq S^{(j)}f^{-1} \\ |r_1| \leq R_1^{(j)}f^{-1}, \, |r_2| \leq R_2^{(j)}f^{-1}\\  sc - r_1r_2  = \square\\  N \mid u^2d_1d_2 c + u(d_1r_2 + d_2r_1) + s}}1.
\end{equation}
We have now completed the transformation of the off-diagonal term; the exponential sums have been bounded by  \eqref{gettingcloser} together with \eqref{case1} if $(j_1, j_2) \not= (2, 2)$ and \eqref{case2} if $(j_1, j_2) = (2, 2)$. It is now clear that a diophantine problem has emerged that we will solve in the next section, at least for such $x \approx u/N$ that cannot be well-approximated. 

\section{Counting congruences}

\begin{lemma} For arbitrary real numbers $C, S, R, \tilde{R} \geq 1$ and positive integers $d_1, d_2, u, N$ let
\begin{displaymath}
  \mathcal{A}_{d_1, d_2}(C, S, R, \tilde{R}; u, N) :=  \{C \leq c < 2C, \, |s| \leq S,  \,  |r_1| \leq R, \, |r_2| \leq \tilde{R} :   N \mid u^2d_1d_2 c + u(d_1r_2 + d_2r_1) + s\}
\end{displaymath}
and 
\begin{displaymath}
  \mathcal{A}^{\square}(C, S, R, \tilde{R}; u, N) := \{(c, s, r_1, r_2) \in  \mathcal{A}_{1, 1}(C, S, R, \tilde{R}; u, N) : sc-r_1 r_2 = \square\}. 
\end{displaymath}
Let $0 < H  \leq N$ be another parameter and assume
\begin{equation}\label{dirichlet}
\left| \frac{u}{N} - \frac{a}{q}\right| \leq \frac{1}{qH}
\end{equation}
for some $q \leq H$ and $(a, q) = 1$. Then
\begin{equation}\label{boundcase1}
 \#\mathcal{A}_{d_1, d_2}(C, S, R, \tilde{R}; u, N)\ll C\min(R, \tilde{R})\left(\frac{S(d_1 \tilde{R} + d_2R)}{N} + \frac{Sq}{N} + \frac{(d_1\tilde{R} + d_2 R)^2}{qH} + \frac{d_1\tilde{R}+d_2R}{q} +1\right)
\end{equation}
and
\begin{equation}\label{boundcase2}
\begin{split}
&  \# \mathcal{A}^{\square}(C, S, R, \tilde{R}; u, N) \ll (CSR\tilde{R} qN)^{\varepsilon} \\& \times \left(\frac{CS(\tilde{R}+R)}{N} + \frac{CSq}{N} + \frac{C(\tilde{R}+R)^2}{qH} + \frac{C(\tilde{R}+R)}{q} + C+ \frac{\sqrt{SC}q\min(R, \tilde{R})}{N}\right).
\end{split}
\end{equation}
\end{lemma}

\textbf{Proof.} An upper bound for $ \#\mathcal{A}_{d_1, d_2}(C, S, R, \tilde{R}; u, N)$ is the number of triplets $(c,r_1, r_2)$ such that
\begin{displaymath}
   \left\| \frac{u}{N}(d_1r_2 + d_2r_1)  -\frac{u^2d_1d_2c}{N} \right\| \leq  \frac{S}{N}.\end{displaymath}
By \eqref{dirichlet}, we can relax this condition to 
\begin{displaymath}
  \left\| \frac{a}{q}(d_1r_2 + d_2r_1) -\frac{u^2d_1d_2c}{N}\right\| \leq \frac{S}{N} +  \frac{d_1\tilde{R}+d_2R}{qH}.
\end{displaymath}
Hence we need to count the number of $C \leq c <2C, |r_1| \leq R, |r_2| \leq \tilde{R}, n \in \Bbb{Z}$ such that
\begin{displaymath}
  a(d_1r_2 + d_2r_1) - qn \in \mathcal{I}(c) := \frac{u^2d_1d_2cq}{N}  + \left[-\frac{Sq}{N} -  \frac{d_1\tilde{R}+d_2R}{H}, \frac{Sq}{N} +  \frac{d_1\tilde{R}+d_2R}{H}\right]. 
\end{displaymath}
The union of the intervals $\mathcal{I}(c)$, $C \leq c < 2C$, covers at most 
\begin{displaymath}
 \ll  C \left(\frac{Sq}{N} +  \frac{d_1\tilde{R}+d_2R}{H}+1\right)
\end{displaymath}
integers.  For a fixed integer $ w \in \mathcal{I}(c)$ and fixed $(a, q) = 1$, the equation $a(d_1r_2 + d_2r_1) - qn = w$ has at most $O(\frac{d_1\tilde{R}+d_2R}{q}+1)$ solutions $n$, and each of them determines at most $O(\min(R, \tilde{R}))$ pairs $(r_1, r_2)$. This completes the proof of \eqref{boundcase1}. 

In order to prove \eqref{boundcase2}, we write $r := r_1+r_2$ and assume first $4sc \not= r^2$. The same argument as above shows that the number of triplets $C \leq c < 2C$, $|s| \leq S$,  $|r| \leq R+\tilde{R}$ satisfying 
\begin{equation}\label{easycong}
  N \mid u^2c + ur+s
\end{equation}
is   
\begin{equation}\label{collect1} 
  \ll C \left(\frac{Sq}{N} +  \frac{\tilde{R}+R}{H}+1\right)\left(\frac{\tilde{R}+R}{q}+1\right).
 \end{equation}
Each of these determines essentially a quadruple $(c, s, r_1, r_2)$: Indeed, the extra condition $sc - r_1r_2  = \square$ 
gives \begin{displaymath}
  4\square - (2r_1-r)^2 = 4sc - r^2 \neq 0,
\end{displaymath}
and so there are at most $(CSR\tilde{R})^{\varepsilon}$ choices for $r_1$ and $r_2$. 

Next we count quadruples $(c, s, r_1, r_2)$ with $4sc = r^2$. Note that in this instance $s>0$ since $c>0$. Writing $c=m c_{2}^2$ and $s=m s_{2}^2$ where $m=(c,s)$, we have $r = \pm 2 m c_2s_2$. Since $N$ is square-free, \eqref{easycong}   gives 
\begin{equation}\label{newcong}
  (N, m) (u c_2 \pm s_2) \equiv 0 \, ( \textnormal{mod }N).
\end{equation}  
We fix $m$ momentarily and use \eqref{dirichlet} as above. Hence the number of pairs $c_2  \asymp \sqrt{C/m}$, $s_2 \leq \sqrt{S/m}$ satisfying \eqref{newcong} is at most the number of $c_2$ satisfying
\begin{displaymath}
  \left\| \frac{a}{q} (N, m)c_2 \right\| \leq \sqrt{\frac{S}{m}} \frac{(N, m)}{N}   + \sqrt{\frac{C}{m}} \frac{(N, m)}{qH} 
\end{displaymath}
or in other words, the number of pairs $(c_2, b)$ with $a(N, m)c_2 \equiv b$ (mod $q$) and
\begin{displaymath}
|b| \leq  \frac{\sqrt{S}q(N,m)}{\sqrt{m}N} + \frac{\sqrt{C}(N, m)}{\sqrt{m}H}.
\end{displaymath}
We write this as 
\begin{displaymath}
  a \frac{(N, m)}{(q, N, m)} c_2 \equiv b' \, \left(\text{mod } \frac{q}{(q, N, m)}\right), \quad\quad |b'| \leq  \frac{\sqrt{S}q(N,m)}{\sqrt{m}N(q, N, m)} + \frac{\sqrt{C}(N, m)}{\sqrt{m}H(q, N, m)} .
\end{displaymath}
Fixing $b'$ determines $c_2$ modulo $q/(q, N, m)$.  If $b'=0$, then $q/(q, N, m)$ divides $c_2$, and each pair $(m, c_2)$ determines at most one  triple $(c, s, r)$ which in turn gives at most $\min(R, \tilde{R})$ quadruples $(c, s, r_1, r_2)$. Hence our count in this case can be bounded by
 \begin{equation}\label{collect2}
\ll \sum_{\substack{0<m c_2^2 \leqslant C \\ c_2 \equiv 0 \, (\textnormal{mod }q/(q, m))}} \min(R,\tilde{R}) \ll  \min(R, \tilde{R})\sum_{ m \leq C} \frac{(q, m)\sqrt{C}}{q\sqrt{m}} \ll q^{\varepsilon} \min(R, \tilde{R}) \frac{C}{q}.
\end{equation}
If $b'\not = 0$,  we find similarly that the number of the quadruples $(c, s, r_1, r_2)$ counted in the present subcase is at most 
\begin{equation}\label{collect3}
\begin{split}
& \ll \sum_{m \leqslant C} \left( \frac{\sqrt{S}q(N, m)}{\sqrt{m}N(q, N, m)} + \frac{\sqrt{C}(N, m)}{\sqrt{m}H(q, N, m)}\right)\left( \sqrt{\frac{C}{m}} \left(\frac{q}{(q, N, m)}\right)^{-1}+1\right) \min(R,\tilde{R}) \\
& \ll \sum_{m \leqslant C} \left(\frac{(N, m)}{m} \left(\frac{\sqrt{SC}}{N} + \frac{C}{Hq}\right) + \frac{(N, m)}{\sqrt{m}} \left(\frac{\sqrt{S}q}{N} + \frac{\sqrt{C}}{H}\right)\right)\min(R, \tilde{R})\\
& \ll (NC)^{\varepsilon}   \left(\frac{\sqrt{SC} q}{N} + \frac{C}{H}\right) \min(R, \tilde{R}). 
\end{split}
\end{equation}
Combining \eqref{collect1} (multiplied by $(CSR\tilde{R})^{\varepsilon}$), \eqref{collect2} and \eqref{collect3} and using $q \leq H$ completes the proof of \eqref{boundcase2}. 

\section{Proof of Theorem 1}

We are now ready to complete the analysis of the off-diagonal terms $\mathcal{Q}_{\kappa}^{\text{off}-\Delta}(z)$ and $\mathcal{Q}(z)$. Lemma 11 below is only a matter of careful book-keeping.  We continue to assume (as we may by Dirichlet's approximation theorem) that
\begin{equation}\label{assume}
  x = \frac{u}{N} + O\left(\frac{1}{N}\right) = \frac{a}{q} + O\left(\frac{1}{N} + \frac{1}{qH}\right)
\end{equation}
for some $(a, q) = 1$, $q \leq H$ and a parameter $H \leq N$ to be chosen later. 
The plan is to substitute \eqref{boundcase1} into \eqref{case1} and \eqref{boundcase2} into \eqref{case2} both of which will be inserted into \eqref{gettingcloser}, the former if $(j_1, j_2) \not= (2, 2)$, and the latter otherwise. We recall the notation \eqref{sizes} and \eqref{sizeS}. Thus, \eqref{case1} can be bounded by
\begin{displaymath}
P^{\varepsilon} t^{\ast} L^j \min\left(\frac{d_1}{d_2}, \frac{d_2}{d_1}\right)^{1/2} \left(\frac{(t^{\ast})^3 L^{3j}(d_1d_2)^{1/2}}{N} + \frac{(t^{\ast})^2 L^{2j}q}{N} + \frac{(t^{\ast})^2 L^{2j} d_1d_2 C}{qH} + \frac{t^{\ast} L^j (d_1d_2)^{1/2} C}{q} + C\right).
\end{displaymath}
To simplify the formula, we assume that 
\begin{equation}\label{boundH}
   t^{\ast} L^2 \leq H \leq N.
\end{equation}   
Recalling the convention $j = (j_1+j_2)/2$, we get by (\ref{boundH}) that $L^j(d_1d_2)^{1/2} \leqslant L^2$ and therefore the previous display can be bounded by (recall $q \leq H$)
\begin{displaymath}
P^{\varepsilon} t^{\ast} L^j \min\left(\frac{d_1}{d_2}, \frac{d_2}{d_1}\right)^{1/2} \left(  \frac{(t^{\ast})^2 L^{2j}H}{N} +  \frac{t^{\ast} L^j (d_1d_2)^{1/2} C}{q} + C\right).
\end{displaymath}
We substitute this into \eqref{gettingcloser} (assuming $(j_1, j_2) \not= (2, 2)$) and obtain
\begin{displaymath}
\begin{split}
 & P^{\varepsilon} (t^{\ast})^2 \frac{Z^{1/2}}{N^{1/2}} \sum_{\substack{i_1, i_2, j_1, j_2 \in \{0, 1, 2\}\\ 1 \leq i_1 + j_1 \leq 2\\ 1 \leq i_2 + j_2 \leq 2\\ j_1+j_2 < 4}}
   \sum_{\substack{d_1 \in \Lambda_{i_1}\\ d_2 \in \Lambda_{i_2}}} \frac{L^{\theta(4 - i_1- j_1 - i_2 -j_2)} }{(d_1d_2)^{3/4}}  L^{(j_1+j_2)/4} \min\left(\frac{d_1}{d_2}, \frac{d_2}{d_1}\right)^{1/2} \\
   & \times \left(  \frac{(t^{\ast})^2 L^{j_1+j_2}H}{N} +  \frac{t^{\ast} L^{3(j_1+j_2)/4} (d_1d_2)^{1/4} Z^{1/2}}{qN^{1/2}} + \frac{Z^{1/2}L^{(j_1+j_2)/4}}{N^{1/2} (d_1d_2)^{1/4}}\right).
\end{split}
\end{displaymath}
At this point we see that the last term in second line produces a total contribution of  $L^2$ if $j_1 = j_2 = 2$,   but we can only win the game if we get a bound $L^{2-\delta}$ which is why we presently exclude this case.  To bound the previous display efficiently requires a case by case analysis according to the possible values of $i_1, i_2, j_1, j_2$, which is completely straightforward. The worst case is $j_1 + j_2 = 3$, in which case $0 \leq i_1 + i_2 \leq 1$; this gives the bound
\begin{equation}\label{finalcase1}
    P^{\varepsilon} (t^{\ast})^2 L^{\theta}  \left(  \frac{(t^{\ast})^2 Z^{1/2}L^{15/4}H}{N^{3/2}} +  \frac{t^{\ast} L^3 Z}{qN} + \frac{Z L^{3/2}}{N}\right).
\end{equation}

Next we assume $j_1 = j_2 = 2$, so $d_1=d_2 = 1$. Using \eqref{boundcase2} and the bounds \eqref{sizes} and \eqref{sizeS} for $j=2$ as well as \eqref{boundH} and $q\leqslant H$, we estimate \eqref{case2} by
\begin{displaymath}
\begin{split}
  P^{\varepsilon} \left(\frac{(t^{\ast})^3 L^6}{N} + \frac{(t^{\ast})^2 L^4 q}{N} + \frac{C(t^{\ast})^2L^4}{qH} + \frac{Ct^{\ast} L^2}{q} + C\right)  \ll  P^{\varepsilon} \left( \frac{(t^{\ast})^2 L^4 H}{N} +   \frac{Ct^{\ast} L^2}{q} + C\right).  
\end{split}
\end{displaymath}
Substituting this into \eqref{gettingcloser}, it is readily seen that the contribution of the term corresponding to $j_1=j_2 = 2$ is dominated by \eqref{finalcase1}.  We summarize the previous computations in the following lemma.

\begin{lemma} Assume \eqref{assume} for some $H$ satisfying \eqref{boundH}. Then 
\begin{displaymath}
  |\mathcal{Q}_{\kappa}^{\text{off}-\Delta}(z)| + |\mathcal{Q}(z)| \ll   P^{\varepsilon} (t^{\ast})^2 L^{\theta}  \left(  \frac{(t^{\ast})^2 Z^{1/2}L^{15/4}H}{N^{3/2}} +  \frac{t^{\ast} L^3 Z}{qN} + \frac{Z L^{3/2}}{N}\right).
\end{displaymath}
\end{lemma}

The scene has now been prepared to prove Theorem 1. Combining Lemma 9 and 11, we obtain
\begin{displaymath}
   |\mathcal{Q}_{\kappa}(z)| + |\mathcal{Q}(z)| \ll   P^{\varepsilon} (t^{\ast})^2 L^{\theta}  \left(  \frac{(t^{\ast})^2 Z^{1/2}L^{15/4}H}{N^{3/2}} +  \frac{t^{\ast} L^3 Z}{qN} + \frac{Z L^{3/2}}{N}\right) + P^{\varepsilon} t^{\ast} L. 
\end{displaymath}
Choosing $B = 3$ and observing \eqref{phitrafo}, we can substitute this bound into \eqref{seethat} and \eqref{seethat2}, and the $\kappa$-sum in \eqref{seethat2} is absolutely convergent. This gives the final estimate
\begin{equation}\label{bound1}
  |f^Z(x+iy)| \ll P^{\varepsilon}(t^{\ast})^{5} L^{\theta/2}\left(\frac{t^{\ast} Z^{1/4}L^{7/8}H^{1/2}}{N^{3/4}}  + \frac{(t^{\ast} LZ )^{1/2}}{(qN)^{1/2}}   + \frac{Z^{1/2}L^{-1/4}}{N^{1/2}} \right) + P^{\varepsilon} (t^{\ast})^{9/2} L^{-1/2},
\end{equation}
still  under the assumptions \eqref{assume} and \eqref{boundH}. In this case, the major arc bound in Lemma 8 gives
\begin{equation}\label{bound2}
 |f^Z(x+iy)| \ll P^{\varepsilon}(t^{\ast})^{3/2} \left(\frac{qZ}{N^{3/2}} + \frac{Z}{H^{3/2}} + \frac{(t^{\ast})^{3/2} q}{Z^{1/2}}\right).
\end{equation}
It remains to optimize the parameters.  We can  assume $Z \geq N^{9/10}$, for otherwise we use  \eqref{trivialZ}. The critical terms are the first and third in \eqref{bound1}, and the second in \eqref{bound2}. In order to match these we use $\theta=7/64$, \eqref{boundZ} and choose
\begin{displaymath}
  H := N^{313/457}(t^{\ast})^{-1803/914} \quad\quad L := N^{64/457}(t^{\ast})^{96/457}.
\end{displaymath}
Under our assumption \eqref{ranget}, this satisfies the requirements \eqref{boundH} and \eqref{rangeL}. Then by \eqref{boundZ}, the bound (\ref{bound1}) becomes
\begin{equation}\label{fin-bound1but1}
  |f^{Z}(x+iy)| \ll P^{\varepsilon} \left((t^{\ast})^{\frac{9979}{1828}} N^{-\frac{25}{914}} + (t^{\ast})^{\frac{11181}{1828}} N^{\frac{71}{914}}q^{-\frac{1}{2}}\right)
\end{equation}
which gives
\begin{equation}\label{fin-bound1}
  |f^{Z}(x+iy)| \ll P^{\varepsilon} \left((t^{\ast})^{\frac{9979}{1828}} N^{-\frac{25}{914}} + (t^{\ast})^{\frac{9979}{1828}} N^{\frac{6158}{75405}}q^{-\frac{1}{2}}\right)
\end{equation}
upon using \eqref{ranget} for the second term.  For (\ref{bound2}), we use $Z \geq N^{9/10}$ and  (\ref{boundZ}) to get 
\begin{equation}\label{fin-bound2}
   |f^{Z}(x+iy)| \ll P^{\varepsilon} \left( (t^{\ast})^{\frac{5}{2}}qN^{-\frac{1}{2}} + (t^{\ast})^{\frac{9979}{1828}} N^{-\frac{25}{914}}  + (t^{\ast})^3 q N^{-\frac{9}{20}}\right). 
\end{equation}
We use \eqref{fin-bound1} if $q \geq q_0$ and \eqref{fin-bound2} if $q < q_0$. We have lots of freedom where to choose the cutting line, for example $q_0 = N^{1/3}$ does the job, and we arrive at
\begin{displaymath}
 |f^{Z}(x+iy)| \ll P^{\varepsilon} (t^{\ast})^{\frac{9979}{1828}} N^{-\frac{25}{914}}.
\end{displaymath} 
This completes the proof of Theorem 1. 

\section{A Theorem of Iwaniec-Sarnak}

In this final section we modify the proof of \eqref{IwSa} in order to show that it holds uniformly up to an $N^{\varepsilon}$ factor for Hecke-Maa{\ss} cusp forms $f$ of arbitrary square-free level $N$ thus giving \eqref{lastsection}:
\begin{prop} Let $f$ be an $L^2$-normalized  Hecke-Maa{\ss} cusp form of square-free level $N$ and let $\varepsilon > 0$.  Then
\begin{displaymath}
  \|f \|_{\infty} \ll_{\varepsilon} N^{\varepsilon} (t^{\ast})^{5/12+\varepsilon}.
\end{displaymath}
\end{prop}
The notation and presentation is kept similar to \cite{IS} so that the reader may easily compare and contrast both works. The starting point is the spectral expansion for an automorphic kernel. Let $h$ be an even function that is holomorphic in $|\Im t|\leq \frac{1}{2}+\delta$ and satisfies $h(t) \ll (|t|+1)^{-2-\delta}$ for some $\delta>0$.  Let $k$ be the Harish-Chandra/Selberg transform of $h$. Then for $z \in \mathbb{H}$ we have
\begin{eqnarray*}
K(z):=\sum_{\gamma \in \Gamma_0(N)/\{\pm 1\}} k(z,\gamma z) & = & \frac{h(i/2)}{\text{vol}(X_0(N))} + \sum_{f \in \mathcal{B}(N, \chi)} h(t_f) |f(z)|^2 \\
& + & \sum_{\substack{\chi_1\chi_2 = \chi\\ b \in \mathcal{B}(\chi_1, \chi_2)}} \frac{1}{4\pi} \int_{\mathbb{R}} h(t) |E_{t, \chi_1, \chi_2, b}(z)|^2 dt
\end{eqnarray*}
where   $k(z,w)=k(u(z,w))$ and 
\begin{equation*}
u(z,w):=\frac{|z-w|^2}{4 \Im(z)\Im(w)}.
\end{equation*}
For our application, we choose 
\begin{equation*}
h(t)=4\pi^2 \frac{\cosh(\frac{\pi t}{2}) \cosh(\frac{\pi T}{2})}{\cosh(\pi t)+\cosh(\pi T)}.
\end{equation*}
This function is nonnegative on $\mathbb{R} \cup i \mathbb{R}$ and satisfies  $h(t) \geq 1$ for  $T\leqslant t \leqslant T+1$. By \cite[Lemma 1.1]{IS}, its Harish-Chandra/Selberg transform has the following properties:
\begin{equation}\label{properties}
\begin{split}
|k(u)| & \leqslant 4T^{1/2} u^{-1/4} (u+1)^{-5/4}, \quad u > 0\\
k(u) & = T+O(1+u^{1/2} T^2), \quad 0 \leq u \leq 1.
\end{split}
\end{equation}

Our goal is to extract geometric information from the spectral expansion by averaging over the actions of Hecke operators on $K(z)$.  For positive integers $n$ co-prime with the level $N$, we define 
\begin{equation*}
R_N (n) = \left\{ \left(\begin{array}{cc}
  a & b\\
  c & d
\end{array}\right) : a, b, c, d \in \mathbb{Z}, ad - bc = n \textnormal{ and } 0\leq c \equiv 0 \, (\textnormal{mod } N)\right\}.
\end{equation*}
Then
\begin{equation}
(T_n K(z)) =  \frac{1}{\sqrt{n}} \sum_{g \in R_N(n)} k(z,gz). \label{geometric}
\end{equation}
On the spectral side\footnote{Note that here we are again using the multiplicative properties of the coefficients of our basis of Eisenstein series.}, we get 
\begin{equation}
\begin{split}
(T_n K(z)) = &  \frac{h(i/2) \sigma(n)}{\text{vol}(X_0(N))\sqrt{n}}  + \sum_{f\in \mathcal{B}(N, \chi)} h(t_f) \lambda_f(n) |f(z)|^2 \\
& + \sum_{\substack{\chi_1\chi_2 = \chi\\ b \in \mathcal{B}(\chi_1, \chi_2)}} \frac{1}{4\pi} \int_{\mathbb{R}} h(t)\lambda_{t, \chi_1, \chi_2, b}(n) |E_{t, \chi_1, \chi_2, b}(z)|^2 dt.
 \label{spectral}
 \end{split}
\end{equation} 
One should notice at this point that the spectral expansion of the automorphic kernel   provides a natural connection between the $L^2$-norms of the Hecke-Maass cusp forms $f$ weighted by their Hecke eigenvalues in (\ref{spectral}) with geometric information in (\ref{geometric}) for each $(n,N)=1$.  
Bounding (\ref{geometric}) and averaging over $(n,N)=1$, along with implementing an amplifier, will produce a bound for any individual form $f$ with $T\leqslant t_f\leq T+1$. By \eqref{atkinlehner} we can assume $y \geq \sqrt{3}/(2 N)$. Given $(n,N)=1$ we wish to obtain a non-trivial bound for the right hand side of \eqref{geometric} by estimating the number of $g\in R_N(n)$ which satisfy $u(z,gz)< \delta$ for some $\delta \geq 0$.  By \eqref{properties}  one sees that $k(u) \ll T$ for all $u\leq n^{-4}$.  Therefore, 
\begin{eqnarray}\label{geometricbound}
\sum_{g \in R_N(n)} k(u(z,gz)) & \ll & T M(z,n,n^{-4}) + T^{1/2} \int_{n^{-4}}^{\infty} \delta^{-1/4} (1+\delta)^{-5/4} d M(\delta)
\end{eqnarray} 
where
\begin{displaymath}
M(\delta)=M(z,n,\delta)=\#\{g \in R_N(n)  \mid  u(z,gz)<\delta\}.\label{countingmeasure}
\end{displaymath}

We proceed as in \cite{IS} by obtaining estimates on $M(\delta)$. We write 
$M(\delta)=M_0(\delta)+M_\star(\delta)$ where $M_0$ counts over the set of upper-triangular matrices
\begin{displaymath}
M_0(\delta)=M_0(z,n,\delta)=\#\{g \in R_N(n)  \mid  u(z,gz)<\delta, \textnormal{ and } c=0<a \} \label{M0}
\end{displaymath}
and $M_\star$ counts over those matrices whose lower left hand entry is strictly positive
\begin{displaymath}
M_\star(\delta)=M_\star(z,n,\delta)=\#\{g \in R_N(n)  \mid  u(z,gz)<\delta \textnormal{ and } c>0 \}. \label{Mstar}
\end{displaymath}
Estimates for $M_0(\delta)$ are straightforward and do not depend on the level $N$ since $c=0$. Therefore, one gets
\begin{equation}
M_0(\delta)\ll_{\varepsilon} n^\varepsilon \left(1+\sqrt{n\delta} y\right) \label{Ubound}
\end{equation}
exactly as in \cite[A.10]{IS}. For the matrices $g\in R_N(n)$ with $c>0$ we have
\begin{equation*}
4n c^2 y^2 u(z,gz) = c^2 |(cz+d)z-az-b|^2
\end{equation*}
or $|(cz+d)(cz-a)+n| = 2 c y \sqrt{un}$ which is the same as,
\begin{equation*}
|n-|cz+d|^2+(cz+d)(2cx+d-a)| = 2 c y \sqrt{un}.
\end{equation*}
For $u<\delta$, the imaginary part of the left hand side above gives
\begin{displaymath}
|2cx+d-a|\leqslant 2 \sqrt{\delta n} \label{imaginarypart}.
\end{displaymath}
Therefore, taking the real part produces
\begin{displaymath}
||cz+d|-\sqrt{n}| \leqslant 4 \sqrt{\delta n} \label{realpart}
\end{displaymath}
so that 
\begin{displaymath}
\bigg|\bigg|\frac{a+d}{2}+icy\bigg|-\sqrt{n}\bigg| \leqslant 5 \sqrt {\delta n}.
\end{displaymath}
Combining the above results and introducing new variables $A=\frac{a-d}{2}$ and $D=\frac{a+d}{2}$ yields the system of equations
\begin{eqnarray}
A & = & cx+O(\sqrt{n\delta})\label{system}, \\ 
\nonumber D^2+c^2y^2 & = & n+O(n(\delta+\sqrt{\delta})),\\ 
\nonumber A^2 & \equiv & D^2-n \, (\textnormal{mod } c), \\ 
\nonumber N & \leqslant & c \ll y^{-1} \sqrt{n(\delta+1)},\\
\nonumber c & \equiv & 0 \, (\textnormal{mod } N).
\end{eqnarray}

For fixed $c$ and $D$, we write $(D^2-n,c)=c_0=c_1^2 c_2$ where $c_1^2$ is the square part of $c_0$ and $c_2$ is square-free.  As remarked in \cite{IS}, there are at most $O(c^{\varepsilon}(1+\frac{c_1 }{c}\sqrt{n\delta}))$ choices of $A$ satisfying (\ref{system}) and then the only remaining matrix entry $b$ is determined.  Therefore, writing $c=mN$ and $W=(Ny)^{-1} \sqrt{n(\delta+1)}$, the number of solutions $A,D,c$ of (\ref{system}) is at most
\begin{equation}
\sum_{1\leqslant m \ll W} (mN)^{\varepsilon} \sum_{c_0|mN} \left(1+\frac{c_1 }{mN}\sqrt{n\delta}\right)\sum_{\substack{D^2+m^2(Ny)^2=n+O(n(\sqrt{\delta}+\delta))\\ (D^2-n,mN)=c_0}}1. \label{totalbound}
\end{equation} 
When $\delta\geq 1$, we have $n + O(n(\sqrt{\delta}+\delta))=O(n\delta)$ and therefore (\ref{totalbound}) above is (note $Ny \geq 1/2$)
\begin{equation}
\begin{split}
& \ll (WN)^{\varepsilon} \left\{\sum_{\substack{m,D \\ D^2+m^2(Ny)^2\ll n\delta}} 1+\sqrt{n\delta}\sum_{1\leqslant m\ll W} \sum_{c_0|mN}\frac{c_1}{mN}\sum_{\substack{(D^2-n)+m^2(Ny)^2\ll n\delta \\ D^2-n\equiv 0 (\textnormal{mod }c_0)}}1 \right\}\\
& \ll (WN)^{\varepsilon} \left\{n\delta+\sqrt{n\delta}\sum_{1\leqslant m \ll W}\sum_{c_0|mN}\frac{\sqrt{n\delta}}{mN}\right\} \ll n\delta (n\delta N)^{\varepsilon}. \label{totalboundlarge}
\end{split}
\end{equation}
When $\delta<1$, then $\sqrt{\delta}$ dominates in $O(n(\sqrt{\delta}+\delta))$ and (\ref{totalbound}) above is
\begin{equation*}
\ll  (WN)^{\varepsilon} \left\{\sum_{\substack{D,m \\D^2+m^2(Ny)^2-n\ll n\sqrt{\delta}}}1 + 
\sqrt{n \delta}\sum_{1\leqslant m \ll W}\sum_{c_0|mN} \frac{c_1}{mN}\sum_{\substack{(D^2-n)+m^2(Ny)^2\ll n\sqrt{\delta}\\ D^2-n\equiv 0 (\textnormal{mod }c_0)}}1\right\}.
\end{equation*} 
By \cite[Lemma 1.4]{IS} we have $\#\{r,s\in \mathbb{Z} \mid |r^2+Ys^2-n|\leqslant n\Delta\} \ll_{\varepsilon} n^{\varepsilon} (n \Delta^{1/2}+1)$ for  $Y \geq 1/2$ and $\Delta\leqslant 1$, hence we can bound the previous display by  
\begin{equation}
\begin{split}
&\ll (WN)^{\varepsilon}\left\{(n\delta^{1/4}+1) + \sqrt{n\delta} \sum_{1\leqslant m\ll W}\sum_{c_0|mN}\frac{c_1}{mN} \left(1+\frac{\sqrt{n}\delta^{1/4}}{c_1} \right)\right\} \\
&\ll (WN)^{\varepsilon} \left\{ (n\delta^{1/4}+1) + \sqrt{n\delta} \sum_{\ell^2m\ll W}\frac{1}{\ell} + n \delta^{3/4} \right\} \ll (nN)^{\varepsilon}(1+n\delta^{1/4}). \label{totalboundsmall}
\end{split}
\end{equation}
Combining the bounds (\ref{totalboundlarge}), (\ref{totalboundsmall}) and (\ref{Ubound}) gives 
\begin{displaymath}
M(\delta) \ll_{\varepsilon} (nN)^{\varepsilon}\left(1 + n(\delta^{1/4}+\delta^{1+\varepsilon}) + \sqrt{n\delta}y\right).\label{Mbound}
\end{displaymath}
In view of (\ref{geometricbound}), one sees that we have arrived at the following bound for the geometric side upon integration by parts:
\begin{equation}
\sum_{g \in R_N(n)} k(u(z,gz))  \ll_{\varepsilon} (n N)^{\varepsilon} \left(T+T^{1/2}n+T^{1/2}n^{1/2}y\right). \label{1.1bound}
\end{equation}
This bound is the same as in \cite[A.11]{IS} with the exception of an additional $N^\varepsilon$ level factor. Next we implement an amplifier.  

Let $\{\alpha_{\ell}\}$ be arbitrary weights supported on $(\ell,N)=1$ and $\ell \leq L$ for some parameter $L$. Let $f_0$ be the $L^2$-normalized cusp form we wish to bound, and let $t_{f_0} \sim T$. Using the multiplication law \eqref{eig2},  we see
\begin{equation*}
\begin{split}
|f_0(z)|^2 \bigg|\sum_{\ell\leqslant L} \alpha_{\ell} \lambda_{f_0}(\ell)\bigg|^2 & \leq \sum_{f\in \mathcal{B}(N, \chi)} h(t_f) |f(z)|^2 \bigg|\sum_{\ell\leqslant L} \alpha_{\ell} \lambda_f(\ell)\bigg|^2 \\
& \leq \sum_{\ell_1,\ell_2\leqslant L} \alpha_{\ell_1} \overline{\alpha_{\ell_2}}\sum_{d|(\ell_1,\ell_2)}\frac{\chi(d) d}{\sqrt{\ell_1 \ell_2}} \sum_{g\in R_N\left(\frac{\ell_1 \ell_2}{d^2}\right)}k(z,gz)
\end{split}
\end{equation*}
where in the last inequality we added the nonnegative contribution of the constant function and the Eisenstein series.  With the resulting bound (\ref{1.1bound}) for the geometric side, the right hand side above is bounded by
\begin{equation*}
(LN)^\varepsilon \sum_{\ell_1,\ell_2\leqslant L } |\alpha_{\ell_1} \alpha_{\ell_2}| \sum_{d|(\ell_1,\ell_2)}\frac{d}{\sqrt{\ell_1\ell_2}}\left(T+T^{1/2}\left(\frac{\ell_1 \ell_2}{d^2}\right)+T^{1/2}\left(\frac{\ell_1\ell_2}{d^2}\right)^{1/2}y\right).
\end{equation*}
By Cauchy-Schwarz, the first term in parentheses  contributes $T(LN)^{\varepsilon} \sum |\alpha_{\ell}|^2$, the second contributes $T^{1/2}L (LN)^{\varepsilon}(\sum |\alpha_{\ell}|)^2$, while the last term contributes  $T^{1/2}y (LN)^{\varepsilon}(\sum |\alpha_{\ell}|)^2$. Note that the last term contribution is actually sharper than listed in \cite{IS} by an $L^{1/2}$ factor. Thus we get the total bound
\begin{eqnarray*}
|f_0(z)|^2 \bigg|\sum_{\ell\leqslant L} \alpha_{\ell} \lambda_{f_0}(\ell)\bigg|^2\ll (LN)^{\varepsilon} \left(T \sum_{\ell\leqslant L} |\alpha_\ell|^2 +  (L+y)T^{1/2}\left(\sum_{\ell\leqslant L}|\alpha_\ell|\right)^2\right).
\end{eqnarray*}
Now we set, similar to \eqref{ampli},
\begin{equation*}
\alpha_\ell := \left\{\begin{array}{ll}
         \lambda_{f_0}(\ell) \bar{\chi}(\ell) & \textnormal{if } \ell=p\leqslant \sqrt{L}\\
         -\bar{\chi}(\ell) & \textnormal{if } \ell=p^2\leqslant L\\
         0 & \textnormal{otherwise}
       \end{array}\right. 
\end{equation*}
for all $p \nmid N$. By \eqref{sizeampli} and \eqref{boundHecke}  this yields
\begin{eqnarray*}
|f_{0}(z)|^2L (LN)^{-\varepsilon}  \ll_{\varepsilon}  (LN)^{\varepsilon}(TL^{1/2} + (L+y)T^{1/2}L),
%
\end{eqnarray*}
and  choosing $L=T^{1/3}$ 
we get that
\begin{equation}\label{fin1}
f_0(z)\ll (t  N)^{\varepsilon} (t^{5/12} + y^{1/2}t^{1/4})
\end{equation}
where $t = t_{f_0} \geq 2$ (without loss of generality). This is uniform in $N$ and 
  completes the proof of Proposition 1 if $ y \leq t^{1/3}$. If $y > t^{1/3}$, we follow the corrigendum to \cite[Lemma A.1]{IS} in \cite{Sa1}. From the Fourier expansion, \eqref{defW},  \eqref{sign} and \eqref{bessel5} we have 
\begin{displaymath}
  |f_0(z)| \ll |\rho_{f_0}(1)| \left(\sum_{n \leq t/y}  |\lambda_{f_0}(n)|  y^{1/2} \cosh(\pi t/2) | K_{it} (2 \pi n y)| + \sum_{n > t/y} \frac{ |\lambda_{f_0}(n)|}{\sqrt{n}} e^{-ny} \right).
\end{displaymath}
The tail is negligible for $y > t^{1/3}$, and we can clearly assume $y \leq t$. By  \eqref{rho1},  Cauchy-Schwarz and \eqref{boundHecke} we find
\begin{equation}\label{laststep}
  |f_0(z)| \ll_{\varepsilon} N^{-1/2} (Nt)^{\varepsilon} t^{1/2} \left(\sum_{n \leq t/y}  |\cosh(\pi t/2)  K_{it} (2 \pi n y)|^2\right)^{1/2} .
 \end{equation}
for any $\varepsilon > 0$. A uniform asymptotic expansion for the Bessel-$K$ function is stated in \cite[(2)]{Ba1} and proved in \cite{Ba2}:
\begin{displaymath}
  K_{it}(tx) = \frac{\pi \sqrt{2}}{t^{1/3}} e^{-\pi t/2} \left(\frac{\zeta}{x^2-1}\right)^{1/4} \left({\rm Ai}(t^{2/3} \zeta) + O\left(\frac{1}{(1+ t^{2/3}|\zeta|)^{1/4}t}\right)\right), \quad t >  1, x > 0
\end{displaymath}
where $\zeta \in \Bbb{R}$ is given by $\frac{2}{3}\zeta^{3/2} = (x^2-1)^{1/2} - {\rm arcsec}(x)$ and ${\rm Ai}(x) \ll (1+|x|)^{-1/4}$ is the Airy-function. Since $\zeta \sim (x-1)$ for $1/2 \leq x \leq 3/2$,  we find by trivial estimates 
\begin{displaymath}
  \cosh(\pi t/2) K_{it}(w) \ll  \min(t^{-1/3}, |w^2-t^2|^{-1/4})
  \end{displaymath}
for $t> 1, w > 0$. Now we estimate the $n$-sum in \eqref{laststep}. The transition range $2 \pi n y = t + O(t^{1/3})$ contributes at most  
\begin{displaymath}
 \ll  \left(\frac{t^{1/3}}{y} + 1\right) t^{-2/3} 
\end{displaymath}
to the $n$-sum, while the contribution of the remaining  $n$ can be bounded by\begin{displaymath}
\begin{split}
& \ll   \sum_{\substack{n \leq t/y\\ |2\pi n y - t| > t^{1/3}}} \frac{1}{|(2 \pi n y)^2 - t^2|^{1/2}} \\
&   \leq \int_0^{\frac{t-t^{1/3}}{2 \pi y}} \frac{ d\xi}{(t^2  - (2\pi y \xi)^2)^{1/2}} + O\left(\frac{1}{t^{2/3}}\right)+  \int_{\frac{t+t^{1/3}}{2 \pi y}}^{t/y} \frac{ d\xi}{((2\pi y \xi)^2 - t^2)^{1/2}}\\
& \ll \frac{1+ \log(ty) }{y} + \frac{1}{t^{2/3}}
 \end{split}
\end{displaymath}
for $t, y > 1$. Collecting the previous estimates we find
\begin{displaymath}
 | f_0(z) | \ll N^{-1/2} \left(\left(\frac{t}{y}\right)^{1/2} + t^{1/6}\right) (Nt y)^{\varepsilon}
\end{displaymath}
which implies Proposition 1 if $t^{1/3}< y \leq t$.

\end{document}